\documentclass[nohyperref]{article}

\usepackage{microtype}
\usepackage{graphicx}
\usepackage{booktabs} %

\usepackage{hyperref}

\usepackage[accepted]{icml2022}

\usepackage{amsmath}
\usepackage{amssymb}
\usepackage{mathtools}
\usepackage{amsthm}

\usepackage[capitalize,noabbrev]{cleveref}

\theoremstyle{plain}

\theoremstyle{definition}

\theoremstyle{remark}

\usepackage[textsize=tiny]{todonotes}

\usepackage{url}
\usepackage{bm}
\usepackage{amsmath}
\usepackage{amsfonts}       %
\usepackage{booktabs}       %

\usepackage{caption}
\usepackage{subcaption}
\usepackage{sidecap}
\usepackage{adjustbox}

\newcommand{\reals}{{\mathbb{R}}}

\icmltitlerunning{Learning a Large Neighborhood Search Algorithm for Mixed Integer Programs}

\begin{document}

\twocolumn[
\icmltitle{Learning a Large Neighborhood Search Algorithm for Mixed Integer Programs}

\icmlsetsymbol{equal}{*}
\begin{icmlauthorlist}
\icmlauthor{Nicolas Sonnerat}{equal,comp}
\icmlauthor{Pengming Wang}{equal,comp}
\icmlauthor{Ira Ktena}{comp}
\icmlauthor{Sergey Bartunov}{dm,charm}
\icmlauthor{Vinod Nair}{dm,gr}
\end{icmlauthorlist}

\icmlEqualContribution

\icmlaffiliation{comp}{DeepMind, London, UK}
\icmlaffiliation{dm}{Work done  at DeepMind}
\icmlaffiliation{gr}{Google Research, Bangalore, India}
\icmlaffiliation{charm}{Charm Therapeutics, London, UK}

\icmlcorrespondingauthor{Nicolas Sonnerat}{sonnerat@deepmind.com}
\icmlcorrespondingauthor{Pengming Wang}{pengming@deepmind.com}

\icmlkeywords{Machine Learning, ICML}

\vskip 0.3in
]

\printAffiliationsAndNotice{}  %

\begin{abstract}
Large Neighborhood Search (LNS) is a combinatorial optimization heuristic that starts with an assignment of values for the variables to be optimized, and iteratively improves it by searching a large neighborhood around the current assignment. In this paper we consider a learning-based LNS approach for mixed integer programs (MIPs). We train a \emph{Neural Diving} model to generate an initial assignment. Formulating the subsequent search steps as a Markov Decision Process, we train a \emph{Neural Neighborhood Selection} policy to select a search neighborhood at each step, which is searched using a MIP solver to find the next assignment. The policy network is trained using imitation learning. We propose a target policy for imitation that is designed to select the neighborhood containing the optimal next assignment amongst all possible choices for the neighborhood of a specified size. Our approach matches or outperforms all the baselines on five diverse real-world MIP datasets with large-scale instances, including two production applications at a large technology company. It achieves $2\times$ to $37.8\times$ better average primal gap than the best baseline on three datasets at large running times.
\end{abstract}

\section{Introduction}
\label{sec:intro}

Large Neighborhood Search (LNS) \cite{shaw1998using,Pisinger2010} is a powerful heuristic for hard combinatorial optimization problems such as Mixed Integer Programs (MIPs)~\cite{Danna2005RINS, rothberg2007evolutionary, berthold2007rens, ghosh2007dins}, Traveling Salesman Problem (TSP)~\cite{smith2017glns}, Vehicle Routing Problem (VRP)~\cite{shaw1998using, hojabri2018large}, and Constraint Programming (CP)~\cite{Perron2004CPLNS, berthold2012large}. Given a problem instance and an initial feasible assignment (i.e., an assignment satisfying all constraints of the problem) of values to the variables of the problem, LNS searches for a better assignment within a neighborhood of the current one at each iteration. Iterations continue until the search budget (e.g., time) is exhausted. The neighborhood is ``large'' in the sense that it contains too many assignments to tractably search with naive enumeration. Large neighborhoods make the search less susceptible to getting stuck in poor local optima.

Two key choices that determine the effectiveness of LNS are 1) the initial assignment, and 2) the search neighborhood at each iteration. A good initial assignment makes good optima more likely to be reached. A good neighborhood selection policy allows faster convergence to good optima. Domain experts design sophisticated heuristics by exploiting problem structure to find an initial feasible assignment, e.g. for MIPs, \cite{fischetti2005feasibility, berthold2007rens} and to define the neighborhood, e.g.~\citet{Pisinger2010, shaw1998using, Danna2005RINS}. 

\begin{figure*}[t]
    \centering
    \includegraphics[width=\textwidth,trim={.9cm 5.25cm 4cm 3.5cm},clip]{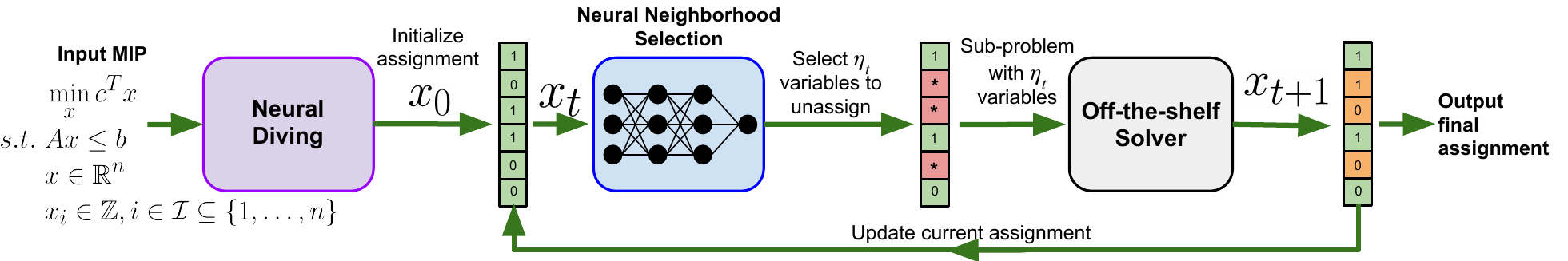}
    \caption{Overview of our approach at test time. The input is a mixed integer program (MIP). \emph{Neural Diving} \cite{nair2020solving} combines a generative model with an off-the-shelf MIP solver to output an initial assignment $x_0$ for the variables $x$ to be optimized.
    At the $t^{\text{th}}$ iteration of LNS the \emph{Neural Neighborhood Selection} policy selects $\eta_t$ variables to be unassigned (indicated by red, with $\eta_t=3$) from the current assignment $x_t$. A sub-MIP defined on those $\eta_t$ variables is solved with a MIP solver to assign them new values (orange) to define the next assignment $x_{t+1}$. Iterations continue until the search budget is exhausted.}
    \label{fig:overview}
\end{figure*}

In this paper we use learned models to make both of these choices. We focus specifically on Mixed Integer Programs to demonstrate the approach, but it can be adapted to other combinatorial optimization problems also. Figure~\ref{fig:overview} gives a summary. To compute an initial feasible assignment of values for the variables, we use \emph{Neural Diving} (section \ref{subsec:background_neural_diving}) proposed in \citet{nair2020solving}, which has been shown to produce high quality assignments quickly. The assignment is computed using a generative model that conditions on the input MIP and defines a distribution over assignments such that ones with better objective values are more probable. To define the search neighborhood at each LNS iteration, we use a \emph{Neural Neighborhood Selection} policy (section \ref{sec:neural_neighborhood_selection}) that, conditioned on the current assignment, selects a subset of the integer variables in the input MIP to unassign their values. The policy's decisions can then be used to derive from the input MIP a smaller ``sub-MIP'' to optimize the unassigned variables. By setting the number of unassigned integer variables sufficiently small, the sub-MIP can be solved quickly using an off-the-shelf solver to compute the assignment for the next LNS step. The policy is trained by imitating an expert neighborhood selection policy (section \ref{subsec:expert_policy}). At each LNS iteration, the expert is designed to select the best assignment in a Hamming ball centered around the current assignment by solving a MIP optimally. The changes in the values of the integer variables between the current and new assignments give the expert's unassignment decisions to use as targets for imitation learning. The expert itself is too computationally expensive for a practical LNS algorithm, but is still tractable for generating imitation training data offline. The neural network trained to approximate it can be orders of magnitude faster, making it practical to use at test time.

Previous works have combined learning with LNS.~\citet{hottung2019neurallns} use an approach complementary to ours for Capacitated VRPs by learning to search the neighborhood, instead of to select it. In our setting, since off-the-shelf MIP solvers can already search neighborhoods, we expect learning to be more useful for neighborhood selection.~\citet{song2020generallns} learn a neighborhood selection policy using imitation learning and reinforcement learning (RL). Their method restricts the neighborhood selection policy to choose fixed, predefined variable subsets, instead of arbitrary subsets in our work. It uses a random neighborhood selection policy to generate training data for imitation learning.~\citet{addanki2020nlns} use RL to learn a policy that unassigns one variable at a time, interleaved with solving a sub-MIP every $\eta$ steps to compute a new assignment. For large MIPs, one neural network policy evaluation per variable can be prohibitively slow. Our approach is scalable -- both selecting an initial assignment and a search neighbourhood at each LNS step are posed as modelling the joint distribution of a large number of simultaneous decisions. This allows us to exploit high-dimensional generative models for scalable training and inference. To demonstrate scalability, we evaluate on real world datasets with large-scale MIPs, unlike earlier works.

\noindent\textbf{Contributions:} 
\begin{enumerate}
    \item We present a scalable learning-based LNS algorithm that combines learned models for computing the initial assignment and for selecting the search neighborhood at each LNS step.
    \item We propose an imitation learning approach to train the neighborhood selection policy using as the imitation target an expert policy that, under certain assumptions, is guaranteed to select the neighborhood containing the optimal next assignment at a given LNS step.
    \item We show results on five diverse large-scale real-world datasets, including two from the production systems of a large technology company. It matches or outperforms all baselines on all of them, with a $2-37.8\times$ improvement over the best baseline with respect to the main performance metric, average primal gap, on three of them.
\end{enumerate}
We have also published source code for the key components of our approach, covering data generation, training, and inference, at \url{https://github.com/deepmind/neural_lns}.

\section{Background}
\label{sec:background}

\subsection{Mixed Integer Programming}
\label{subsec:background_mip}
A Mixed Integer Program is defined as $\min_{x}\{f(x) = c^Tx~~ |~~ Ax \leq b,~~ x_i \in \mathbb{Z}, ~~i \in \mathcal{I}\}$, where $x \in \mathbb{R}^n$ are the variables to be optimized, $A \in \mathbb{R}^{m \times n}$ and $b \in \mathbb{R}^m$ specify $m$ linear constraints, $c \in \mathbb{R}^n$ specifies the linear objective function, and $\mathcal{I} \subseteq \{1,\ldots,n\}$ is the index set of integer variables. If $\mathcal{I} = \emptyset$, the resulting continuous optimization problem is called a \emph{linear program}, which is solvable in polynomial time. A \emph{feasible} assignment is a point $x \in \mathbb{R}^n$ that satisfies all the constraints. A \emph{complete solver} tries to produce a feasible assignment and a lower bound on the optimal objective value, and given sufficient compute resources will find the optimal assignment or prove that there exists no feasible one. A \emph{primal heuristic} (see, e.g., \citealt{berthold2006thesis}) only attempts to find a feasible assignment. This work focuses on primal heuristics as production applications often only require finding a good feasible assignment quickly.

\subsection{Neural Diving}
\label{subsec:background_neural_diving}
Neural Diving \cite{nair2020solving} is a learning-based primal heuristic. It learns a probability distribution for assignments of integer variables of the input MIP $M$ such that assignments with better objective values have higher probability. Assuming minimization, an \emph{energy function} is defined over the integer variables of the problem $x_\mathcal{I} = \{x_i | i \in \mathcal{I}\}$ as 
\begin{equation}
    E(x_{\mathcal{I}}; M) = \begin{cases}
        \hat{f}(x_\mathcal{I}) & \text{if $x_\mathcal{I}$ is feasible,}\\
        \infty & \text{otherwise},\\
    \end{cases}
\end{equation}
where $\hat{f}(x_\mathcal{I})$ is the objective value obtained by substituting $x_\mathcal{I}$ for the integer variables in $M$ and assigning the continuous variables to the solution of the resulting linear program. The target distribution is defined as $p(x_{\mathcal{I}}|M) = \exp(-E(x_{\mathcal{I}};M))/Z(M)$ where $Z(M) = \sum_{x_{\mathcal{I}}'}\exp(-E(x_{\mathcal{I}}';M))$ is the \emph{partition function}. The model is trained to minimize a weighted negative log likelihood loss on the training set $\{(M^{(j)},x_{\mathcal{I}}^{(j)}))\}_{j=1}^N$ of $N$ MIPs and corresponding assignments collected with an off-the-shelf solver. Assignments with better objective values are given larger weights so that the model better approximates the target distribution.

\citet{nair2020solving} represent the MIP as a bipartite graph (see, e.g., \citet{gasse2019exact}) and use a Graph Convolutional Network \cite{battaglia2018relational,gori2005new,scarselli2008graph,hamilton2017inductive,kipf2016semi} to parameterize a conditionally independent model of the joint distribution over integer variable assignments.

Given a MIP at test time, the model's predicted distribution is used to generate multiple \emph{partial} assignments for the integer variables. For each such partial assignment, substituting the values of the assigned variables in $M$ defines a sub-MIP with only the unassigned variables, which is then solved by an off-the-shelf MIP solver to complete the assignment. Neural Diving outputs the best assignment among all such complete assignments. Since each partial assignment can be completed in parallel, Neural Diving is well-suited to exploit parallel computation for faster runtimes. Results in~\citet{nair2020solving} show that Neural Diving is effective in producing high quality assignments quickly on several datasets. See that paper for further details.
\section{Neural Neighborhood Selection}
\label{sec:neural_neighborhood_selection}

\subsection{MDP Formulation}
\label{subsec:mdp}
We consider a contextual Markov Decision Process \cite{abbasi2014online, hallak2015contextual} $\mathcal{M}_z$ parameterized with respect to a \emph{context} $z$, where the state space, action space, reward function, and the environment all depend on $z$. Here we define $z$ to be the parameters of the input MIP, i.e., $z = M = \{A, b, c\}$. The \emph{state} $s_t$ at the $t^{th}$ step of an episode is the current assignment $x_t$ of values for all the integer variables in $M$. The \emph{action} $a_t \in \{0,1\}^{|\mathcal{I}|}$ at step $t$ is the choice of the set of integer variables to be unassigned, specified by one indicator variable per integer variable in $M$ where 1 means unassigned. All continuous variables are labelled as unassigned at every LNS step. For real-world applications the number of integer variables $|\mathcal{I}|$ is typically large ($10^3-10^6$), so the actions are high-dimensional binary vectors. The \emph{policy} $\pi_{\theta}(a_t | s_t, M)$ defines the distribution over actions, parameterized by $\theta$. We use a conditional generative model to represent this high-dimensional distribution over binary vectors (section \ref{subsec:policy_network_architecture}). 

Given $s_t$ and $a_t$, the \emph{environment} derives a sub-MIP $M'_t = \{A'_t, b'_t, c'_t\}$ from $M$ containing only the unassigned integer variables and all continuous variables, and optimizes it. $M'_t$ is computed by substituting the values in $x_t$ of the assigned variables into $M$ to derive constraints and objective function with respect to the rest of the variables. $M'_t$ is guaranteed to have a non-empty feasible set -- the values in $x_t$ of the unassigned variables itself is a feasible assignment for $M'_t$. The set of feasible assignments for $M'_t$ is the search neighborhood for step $t$. The environment calls an off-the-shelf MIP solver, in our case the state-of-the-art non-commercial MIP solver SCIP 7.0.1 \cite{gamrath2020scip}, to search this neighborhood. The output of the solve is then combined with the values of the already assigned variables to construct a new feasible assignment $x_{t+1}$ for $M$. If the solver outputs an optimal assignment for the sub-MIP, then $c^Tx_{t+1} \leq c^Tx_{t}$. The per-step \emph{reward} can be defined using a metric that measures progress towards an optimal assignment \cite{addanki2020nlns}, such as the negative of the \emph{primal gap} \cite{berthold2006thesis} (see equation \ref{eqn:primal_gap}) which is normalized to be numerically comparable across MIPs (unlike, e.g., the raw objective values).

An episode begins with an input MIP $M$ and an initial feasible assignment $x_0$. It proceeds by running the above MDP to perform large neighborhood search until the search budget (e.g., time) is exhausted.

The size of the search neighborhood at the $t^{th}$ step typically increases exponentially with the number of unassigned integer variables $\eta_t$. Larger neighborhoods can make LNS less susceptible to getting stuck at local optima, but it can also be computationally more expensive to search. We treat $\eta_t$ as a hyperparameter to control this tradeoff.

\subsection{Expert Neighborhood Selection Policy}
\label{subsec:expert_policy}
We propose an expert policy that aims to compute the unassignment decisions $a^*_t$ for finding the optimal next assignment $x^*_{t+1}$ across all possible search neighborhoods around $x_t$ given by unassigning any $\eta_t$ integer variables. It uses \emph{local branching} \cite{fischetti03} to compute the optimal next assignment $x^*_{t+1}$ within a given Hamming ball of radius $\eta_t$ centered around the current assignment $x_t$. The minimal set of unassignment decisions $a^*_t$ is then derived by comparing the values of the integer variables between $x_t$ and $x^*_{t+1}$ and labelling only those with different values as unassigned. If the policy $\pi_{\theta}(a_t | x_t, M)$ takes the action $a^*_t$ and the corresponding sub-MIP $M'_t = \{A'_t, b'_t, c'_t\}$ is solved optimally by the environment, then the next assignment will be $x^*_{t+1}$.

Local branching adds a constraint to the input MIP $M$ such that only those assignments for $x_{t+1}$ within a Hamming ball around $x_t$ are feasible. If all integer variables in $M$ are binary, the constraint is:
\begin{equation}
    \sum_{i\in\mathcal{I}: x_{t,i}=0}x_{t+1,i} + \sum_{i\in\mathcal{I}: x_{t,i}=1}(1-x_{t+1,i}) \leq \eta_t, 
\end{equation}
where $x_{t,i}$ denotes the $i^{th}$ dimension of $x_t$ and $\eta_t$ is the desired Hamming radius. The case of general integers can also be handled (see, e.g., slide 23 of \citet{lodi2003LocalBranchingTutorial}). The optimal solution of the MIP with the extra constraint will differ from $x_t$ only on at most $\eta_t$ dimensions, so it is the best assignment across all search neighborhoods for the desired number of unassigned integer variables.

The expert itself is too slow to be directly useful for solving MIPs, especially when the number of variables and constraints are large. Instead it is used to generate episode trajectories from a training set of MIPs for imitation learning. As a one-time offline computation, the compute budget for data generation can be much higher than that of solving a MIP, which enables the use of a slow expert.

\subsection{Policy Network}
\label{subsec:policy_network_architecture}

\noindent\textbf{MIP representation:} Following \citet{nair2020solving} and earlier works (e.g., \citet{gasse2019exact}), we use a bipartite graph representation of a MIP for both Neural Diving and Neural Neighborhood Selection. Variables form one set of nodes in the graph and constraints form the other set. A variable appearing in a constraint is indicated by an edge between the two corresponding nodes. Coefficients in $A$, $b$, and $c$ are encoded as features of the corresponding edges, constraint nodes, and variable nodes, respectively. Additional features can be included (e.g., the linear relaxation solution as variable node features) -- we use the feature set proposed in \citet{gasse2019exact}. For the Neural Neighborhood Selection policy network, we additionally use a fixed-size window of past variable assignments as variable node features. The window size is set to 3 in our experiments. 

\noindent\textbf{Network architecture:} We use a Graph Convolutional Network to represent the policy. Let the input to the GCN be a graph $G = (\mathcal{V}, \mathcal{E}, \mathcal{A})$ defined by the set of nodes $\mathcal{V}$, the set of edges $\mathcal{E}$, and the graph adjacency matrix $\mathcal{A}$. In the case of MIP bipartite graphs, $\mathcal{V}$ is the union of $n$ variable nodes and $m$ constraint nodes, of size $K := |\mathcal{V}| = n+m$. $\mathcal{A}$ is an $K\times K$ binary matrix with $\mathcal{A}_{ij}=1$ if nodes indexed by $i$ and $j$ are connected by an edge, 0 otherwise, and $\mathcal{A}_{ii}=1$ for all i. Let $U \in \mathbb{R}^{K \times D}$ be the matrix containing $D$-dimensional feature vectors of all nodes as rows.

An $L$-layer GCN is defined as follows:
\begin{align}
    Z^{(0)} &= U\\
    Z^{(l+1)} &= \mathcal{A}g_{\phi(l)}(Z^{(l)}), \quad l=0, \ldots, L-1,
\end{align}
where $Z^{(l)} \in \mathbb{R}^{K \times H^{(l)}}$ is the matrix of $H^{(l)}$-dimensional \emph{node embeddings} for the $K$ nodes as rows, and $g_{\phi(l)}()$ is a Multi-Layer Perceptron (MLP) \citep{Goodfellow-et-al-2016} with learnable parameters $\phi(l) \in \theta$ for the $l^\text{th}$ layer applied row-wise to $Z^{(l)}$. The $L^\text{th}$ layer's node embeddings can be used as input to another MLP that computes the outputs for the final prediction task. %

The policy is a conditionally independent model
\begin{equation}
    \pi_{\theta}(a_t | x_t, M) = \prod_{i \in \mathcal{I}} p_{\theta}(a_{t,i} | x_t, M),
\end{equation}
which predicts the probability of $a_{t,i}$, the $i^\text{th}$ dimension of $a_t$, independently of its other dimensions conditioned on $M$ and $x_t$ using the Bernoulli distribution $p_{\theta}(a_{t,i} | x_t, M)$. Its success probability $\mu_{t,i}$ is computed as 
\begin{align}
    \lambda_{t,i} &= \text{MLP}(v_{t,i}; \theta),\\
    \mu_{t,i} &= p_{\theta}(a_{t,i} = 1 | x_t, M) = \frac{1}{1+\exp(-\lambda_{t,i})},
\end{align}
where $v_{t,i} \in \reals^{H^((L)}$ is the $L^\text{th}$ layer embedding computed by a GCN for the node corresponding to $x_{t,i}$, and $\lambda_{t,i} \in \reals$.

\subsection{Training}
\label{subsec:training}
Given a training set $\mathcal{D}_{\text{train}} = \{(M^{(j)},x_{1:T_j}^{(j)},a_{1:T_j}^{(j)}))\}_{j=1}^N$ of $N$ MIPs and corresponding expert trajectories, the model parameters $\theta$ are learned by minimizing the negative log likelihood of the expert unassignment decisions:
\begin{equation}\label{eq:maxlikelihood}
    L(\theta) = -\sum_{j=1}^N \sum_{t=1}^{T_j} \log \pi_{\theta}(a_{t}^{(j)} | x_{t}^{(j)}, M^{(j)}),
\end{equation}
where $M^{(j)}$ is the $j^{th}$ training MIP instance, $\{x_{t}^{(j)}\}_{t=1}^{T_j}$ are the feasible assignments for the variables in $M^{(j)}$, and $\{a_{t}^{(j)}\}_{t=1}^{T_j}$ are the corresponding unassignment decisions by the expert in a trajectory of ${T_j}$ steps.

\subsection{Using the Trained Model}
\label{subsec:using_trained_model}
Given an input MIP, first Neural Diving is applied to it to compute the initial feasible assignment. An episode then proceeds as described in section \ref{subsec:mdp}, with actions sampled from the trained model.

\noindent\textbf{Sampling actions:} Directly sampling unassignment decisions from the Bernoulli distributions output by the model often results in sets of unassigned variables that are much smaller than a desired neighborhood size. This is due to highly unbalanced data produced by the expert (typically most of the variables remain assigned), which causes the model to predict a low probability of unassigning each variable. Instead we construct the unassigned variable set $U$ sequentially, starting with an empty set, and at each step adding to it an integer variable $x_{t,i}$ with probability proportional to $(p_{\theta}(a_{t,i} = 1 | x_t, M) + \epsilon)^{\frac{1}{\tau}} \cdot \mathbb{I}[x_{t,i} \notin U]$ with $\epsilon > 0$ to assign nonzero selection probability for all variables. Here, $\tau$ is a temperature parameter. This ensures that $U$ contains exactly the desired number of unassigned variables.

\noindent\textbf{Adaptive neighborhood size:} The number of variables unassigned at each step is chosen in an adaptive manner \cite{stuckey2021coursera}. The initial number is set as a fraction of the number of integer variables in the input MIP. At a given LNS step, if the sub-MIP solve outputs a provably optimal assignment, the fraction for the next step is increased by a factor $\alpha > 1$. If the sub-MIP solve times out without finding a provably optimal assignment, the fraction for the next step is divided by $\alpha$. This allows LNS to adapt the neighborhood size according to difficulty of the sub-MIP solves.

\section{Evaluation Setup}
\label{sec:evaluation}

\subsection{Datasets}
\label{subsec:datasets}
We evaluate our approach on five datasets: Neural Network Verification, Electric Grid Optimization, Production Packing, Production Planning, and MIPLIB. The first four are homogeneous datasets in which the instances are from a single application, while MIPLIB \cite{gleixner2019miplib} is a heterogeneous public benchmark with instances from many, often unrelated, applications. They contain large-scale MIPs with thousands to millions of variables and constraints (figure~\ref{fig:presolved_mip_sizes} in Technical Appendix). In particular, Production Packing and Planning datasets are obtained from a large technology company's production systems. See Technical Appendix section \ref{subsec:appendix_datasets} for details. All five datasets were split into training, validation, and test sets, each consisting of 70\%, 15\%, and 15\% of total instances, respectively. We train a separate model on each dataset, and evaluate it on the corresponding test set's MIPs. We also report in section \ref{subsec:neurips21_results} preliminary results for datasets from the NeurIPS'21 Machine Learning for Combinatorial Optimization competition.

\subsection{Metrics}
\label{subsec:metrics}
We follow the evaluation protocol of~\citet{nair2020solving} and report two metrics, the \emph{primal gap} and the fraction of test instances with the primal gap below a threshold, both as a function of time. The \emph{primal gap} $\gamma(t)$ at time $t$ is the normalized difference between the objective value achieved by an algorithm under evaluation at $t$ and a precomputed best known objective value $f(x^*)$~\citep{berthold2006thesis}:
\begin{equation}
    \label{eqn:primal_gap}
    \gamma(t) = 
    \begin{cases}
        1, & \text{if } f(x_t) \cdot f(x^*) < 0  \\
         & \text{or no solution at time } $t$, \\
        \frac{|f(x_t) - f(x^*)|}{\max\{ |f(x_t)|, |f(x^*)| \}}, & \text{otherwise}.
    \end{cases}
\end{equation}
We average primal gaps over all test instances at a given time and refer to this as \emph{average primal gap}, and plot it as a function of running time.

Applications typically specify a threshold on the gap between an assignment's objective value and a lower bound, below which the assignment is deemed close enough to optimal to stop the solve. The dataset-specific gap thresholds are given in table \ref{tab:target_opt_gaps}. We apply these thresholds to the primal gap to decide when a MIP is considered solved. We plot the fraction of solved test instances as a function of running time, which we refer to as a \emph{survival curve}.

As in \citet{nair2020solving}, we use \emph{calibrated time} to measure running time. It reduces the variance of time measurements on a shared compute cluster needed for a large-scale evaluation. See the Technical Appendix for details.

\subsection{Baselines}
\label{subsec:baselines}
We compare our approach to three baselines:
\begin{enumerate}
    \item Random Neighborhood Selection (RNS), where the integer variables to unassign are selected uniformly randomly (referred to as the \emph{Random Destroy method} in \citet{Pisinger2010}), with an adaptive neighbourhood size as explained in section \ref{subsec:using_trained_model}. We use Neural Diving to initialize the feasible assignment.
    \item Neural Diving, as described in \citet{nair2020solving} and section \ref{subsec:background_neural_diving}, which achieved the previous best primal gap results on the datasets in section \ref{subsec:datasets}.
    \item SCIP 7.0.1 with its hyperparameters (``metaparameters" for the presolve, cuts, and heuristics components) tuned for each dataset separately using grid search to achieve the best validation set average primal gap curves. SCIP is a complete solver that uses state-of-the-art primal heuristics. By tuning SCIP's hyperparameters to minimize average primal gap quickly, we aim to make SCIP behave more like a primal heuristic.
\end{enumerate}

\noindent\textbf{Use of parallel computation:} Neural Diving can naturally exploit parallel computation for faster performance. This advantage carries over to Neural Neighborhood Selection as well when combined with Neural Diving by using parallel LNS runs initialized with multiple feasible assignments. We evaluate Neural Diving and combinations of Neural Diving with Random or Neural Neighborhood Selection in the parallel setting. All of these primal heuristics are given the same amount of parallel compute resources in experiments. SCIP is evaluated only in the single core setting, as the main focus of this work is to evaluate the benefit of easily parallelizable primal heuristics.

\section{Results}
\label{sec:results}

Figure \ref{fig:primal_gap} shows that on all five datasets, combining Neural Diving and Neural Neighbourhood Selection (ND + NNS) significantly outperforms SCIP on the test instances, in some cases substantially. On Production Packing, the final average primal gap is almost two orders of magnitude smaller, while on Neural Network Verification and Production Planning it is more than $10\times$ smaller. On all datasets except MIPLIB, the advantage of ND + NNS over SCIP is substantial even at smaller running times.

ND + NNS outperforms Neural Diving alone on all datasets, with $10-100\times$ smaller gap on Production Packing, Electric Grid Optimization, and Neural Network Verification. Neural Diving quickly reduces the average primal gap early on, but plateaus at larger running times. ND + NNS overcomes this limitation, reducing the average primal gap significantly with more running time. On MIPLIB, Neural Diving shows a better gap curve initially, before being overtaken by ND + NNS after about $10^3$ seconds.

Combining Neural Diving with Random Neighbourbood Selection (ND + RNS) is a strong baseline on all datasets except Electric Grid Optimization. It is only slightly worse than ND + NNS on Neural Network Verification and MIPLIB. But on Production Planning, Production Packing, and Electric Grid Optimization, ND + NNS achieves a final average primal gap that is smaller by roughly $2.0\times$, $13.9\times$, and $37.8\times$, respectively. Note that ND + RNS is not better than Neural Diving alone on all datasets, but ND + NNS is.

\subsection{Survival Curves}
\label{subsec:survival_plots}
\begin{figure*}
    \centering
    \begin{subfigure}[t]{0.3\textwidth}
        \centering
        \includegraphics[width=\linewidth]{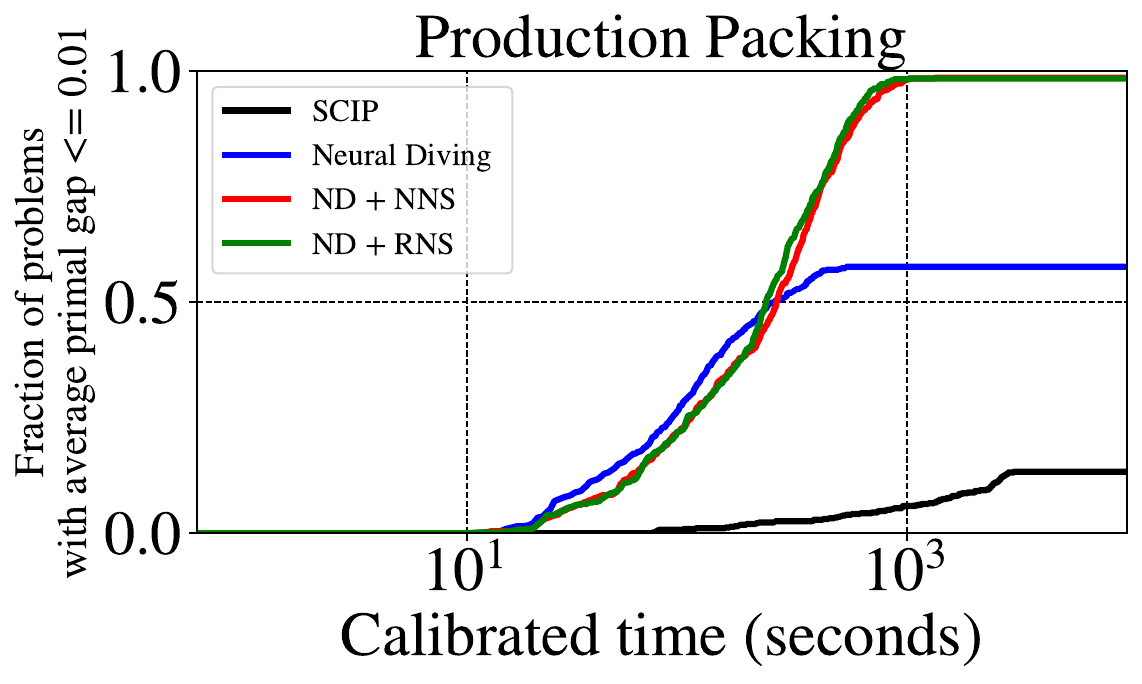}
        \includegraphics[width=\linewidth]{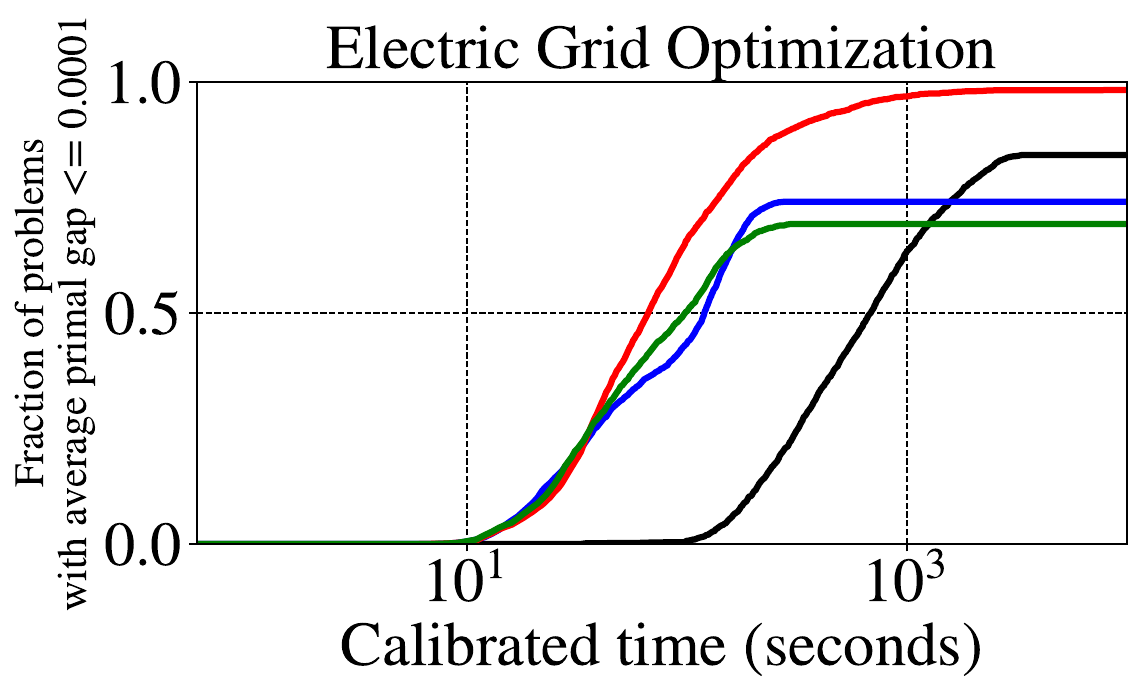}
        \includegraphics[width=\linewidth]{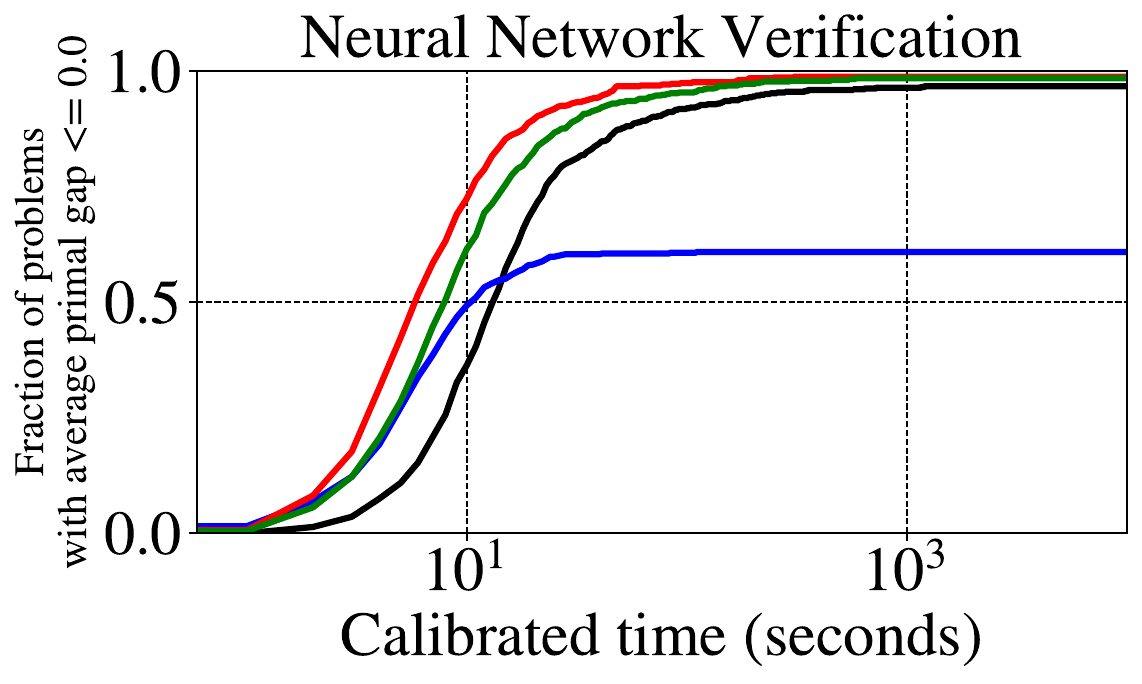}
        \includegraphics[width=\linewidth]{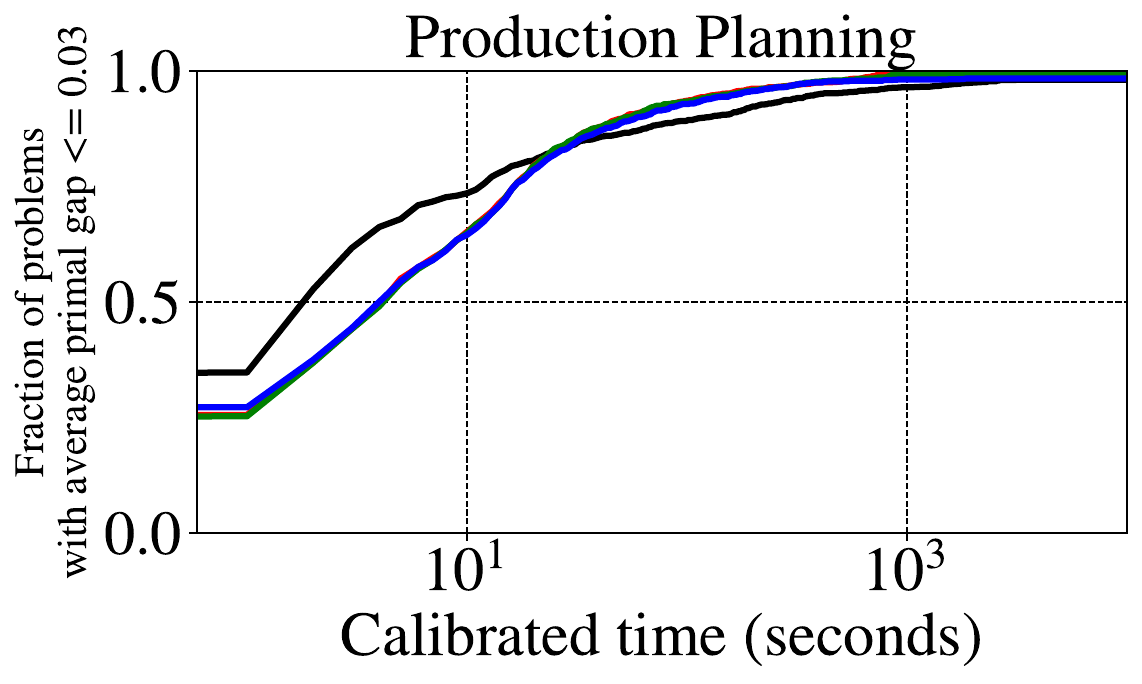}        
        \includegraphics[width=\linewidth]{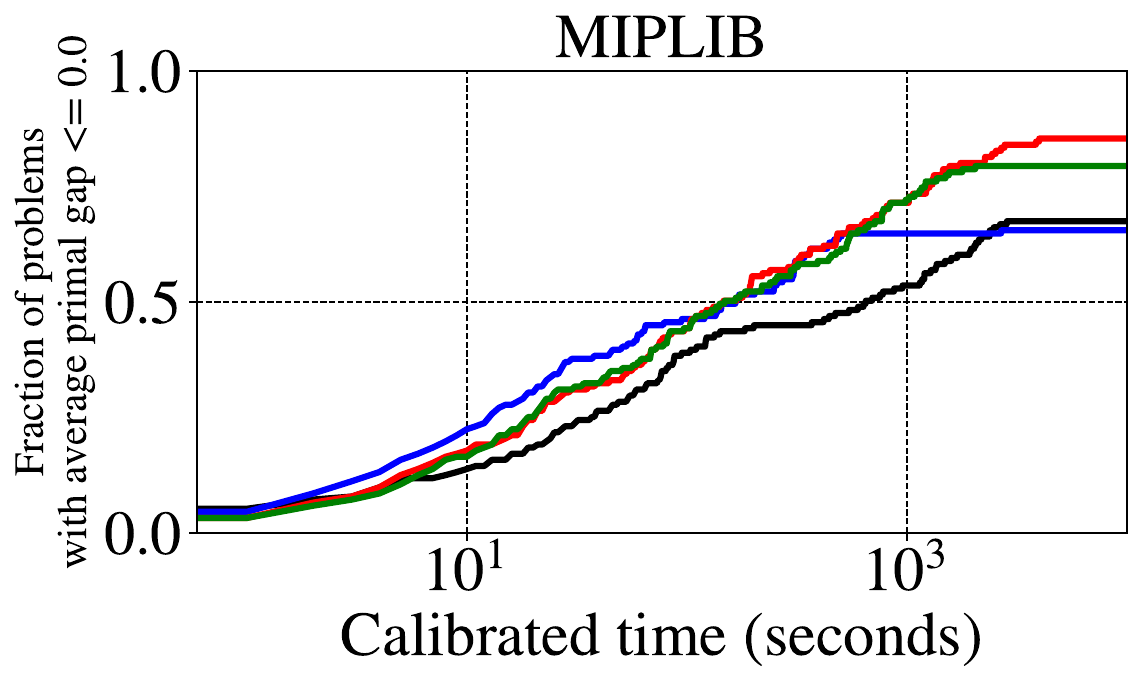} 
        \caption{Fraction of test set instances with primal gap below a dataset-specific threshold, as a function of running time for five datasets. (Note: For Production Planning, several curves closely overlap.)}
        \label{fig:survival}
    \end{subfigure}
    \hfill
    \begin{subfigure}[t]{0.3\textwidth}
        \centering
        \includegraphics[width=\linewidth]{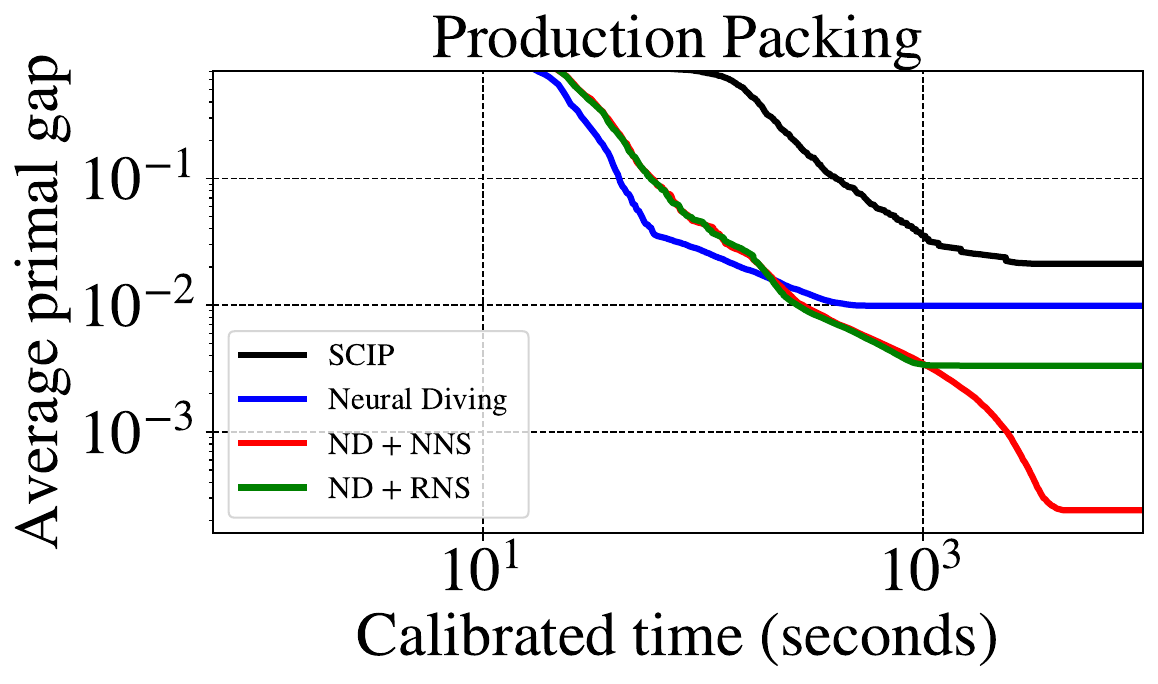} 
        \includegraphics[width=\linewidth]{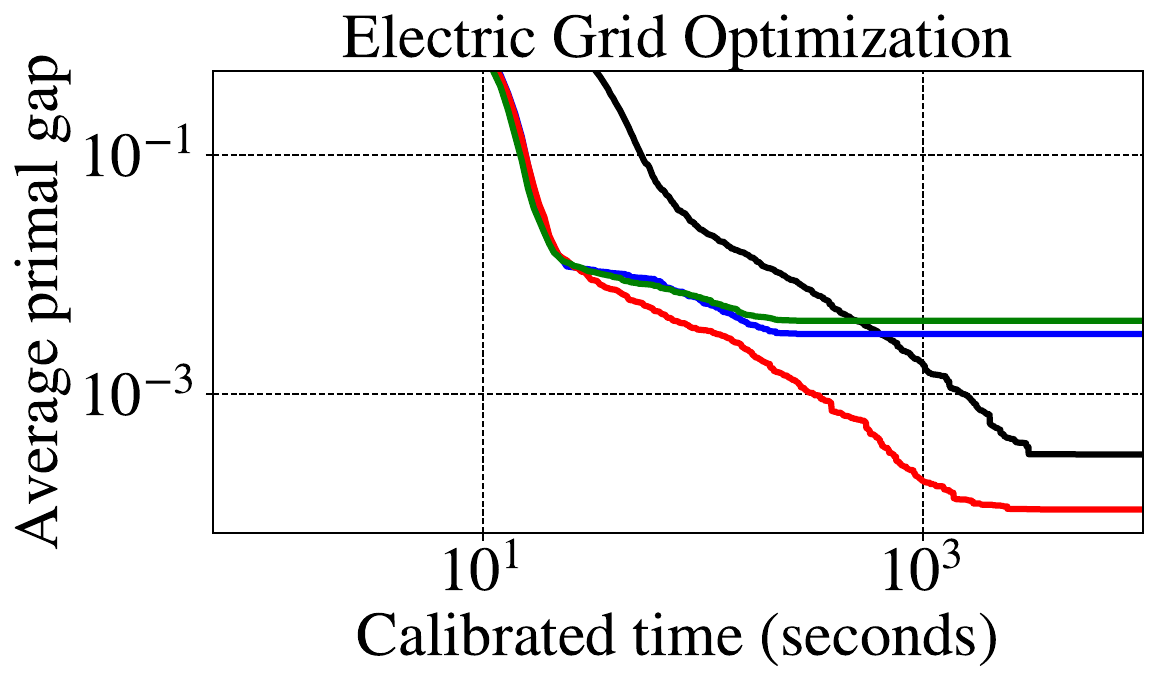} 
        \includegraphics[width=\linewidth]{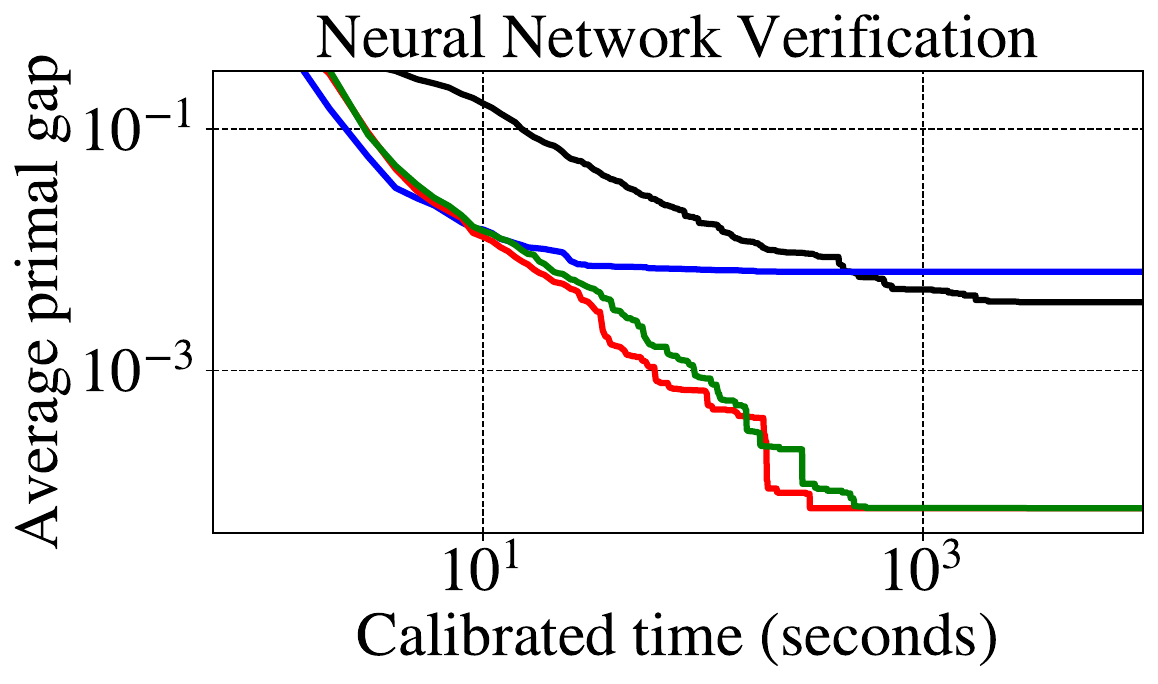}
        \includegraphics[width=\linewidth]{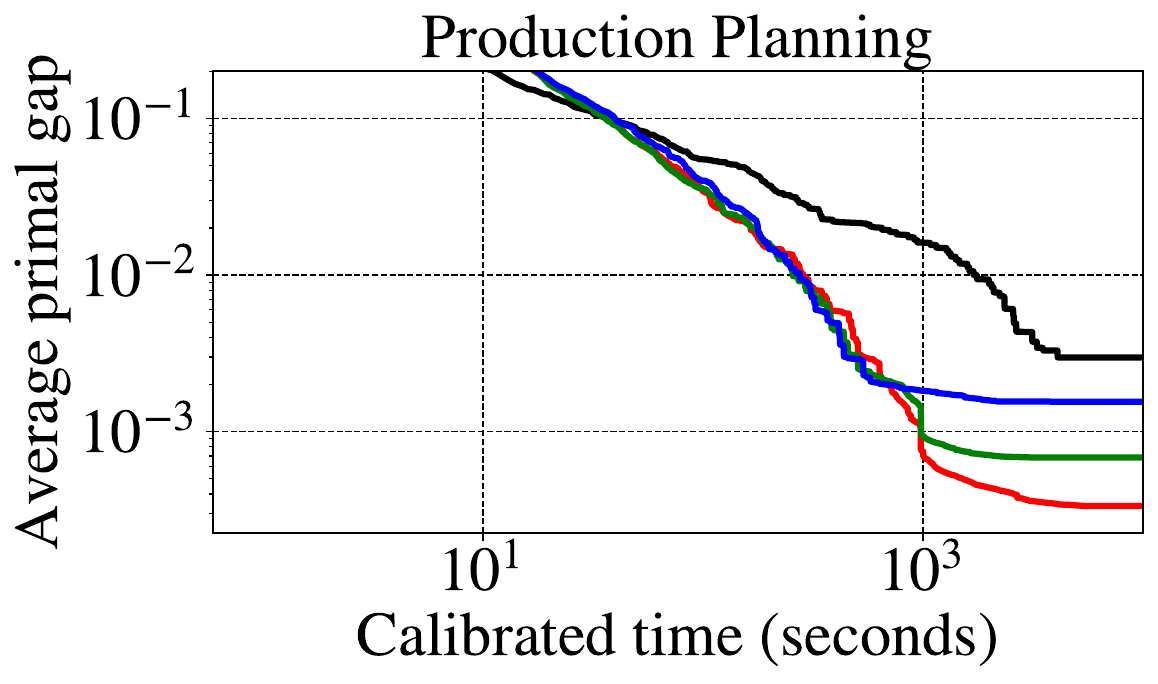}         
        \includegraphics[width=\linewidth]{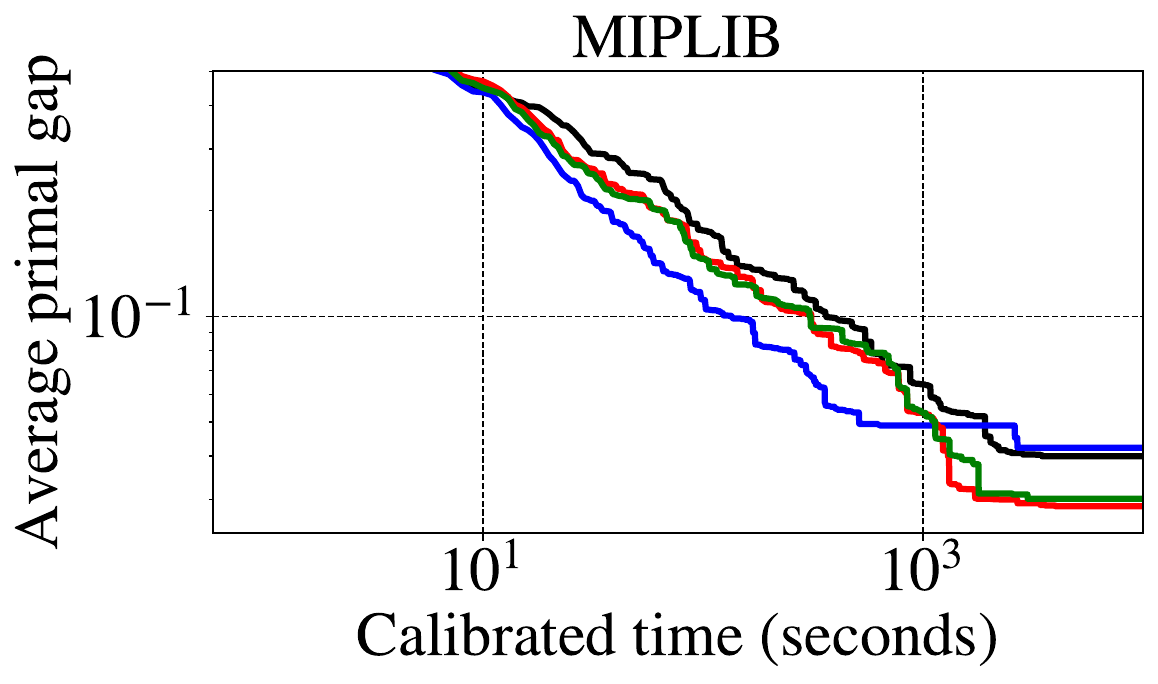} 
        \caption{Test set average primal gap (see section \ref{subsec:metrics}, lower is better) as a function of running time for five datasets.}
        \label{fig:primal_gap}
    \end{subfigure}
    \hfill
    \begin{subfigure}[t]{0.3\textwidth}
        \centering
        \includegraphics[width=\linewidth]{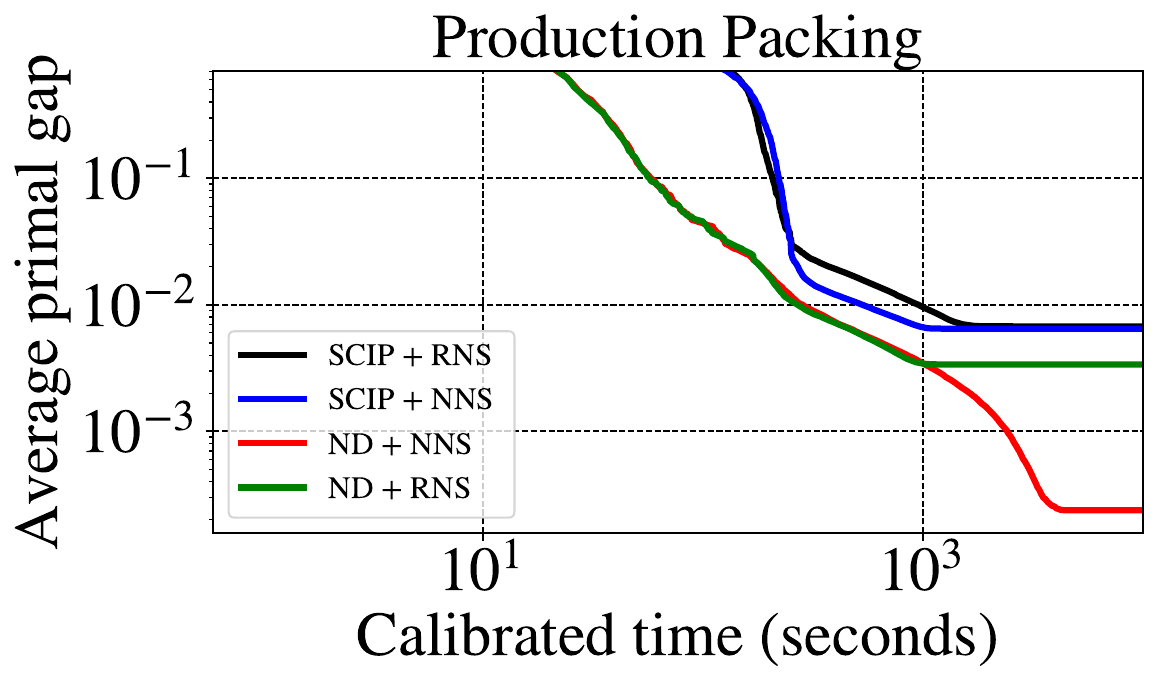}
        \includegraphics[width=\linewidth]{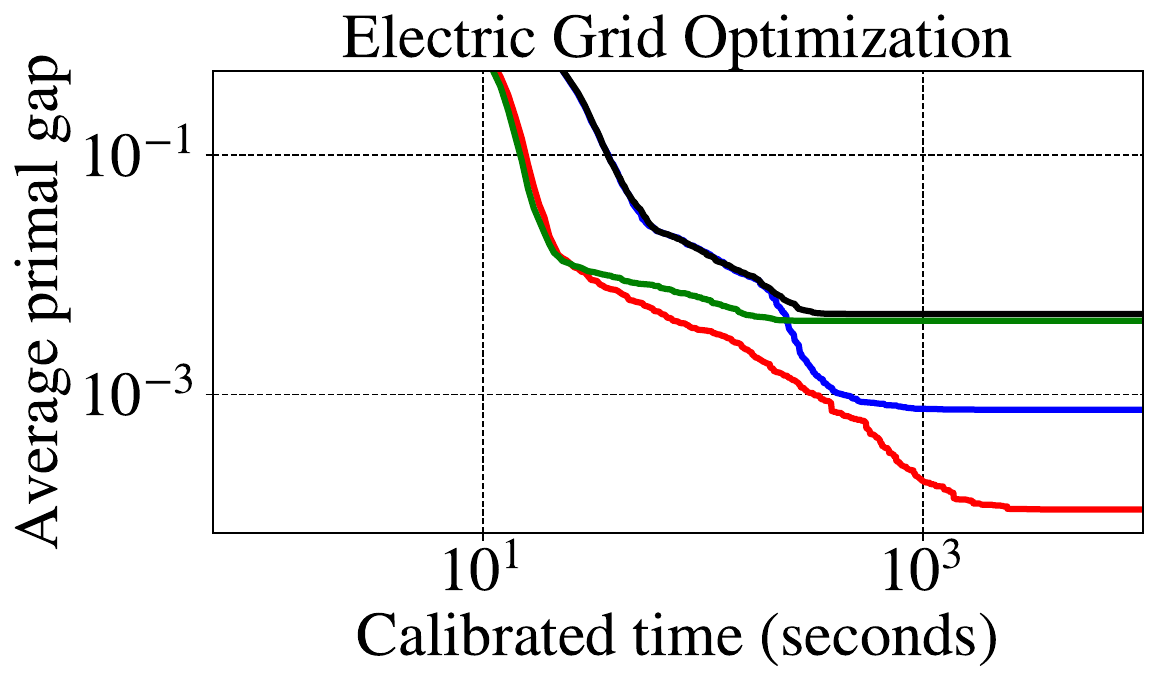}
        \includegraphics[width=\linewidth]{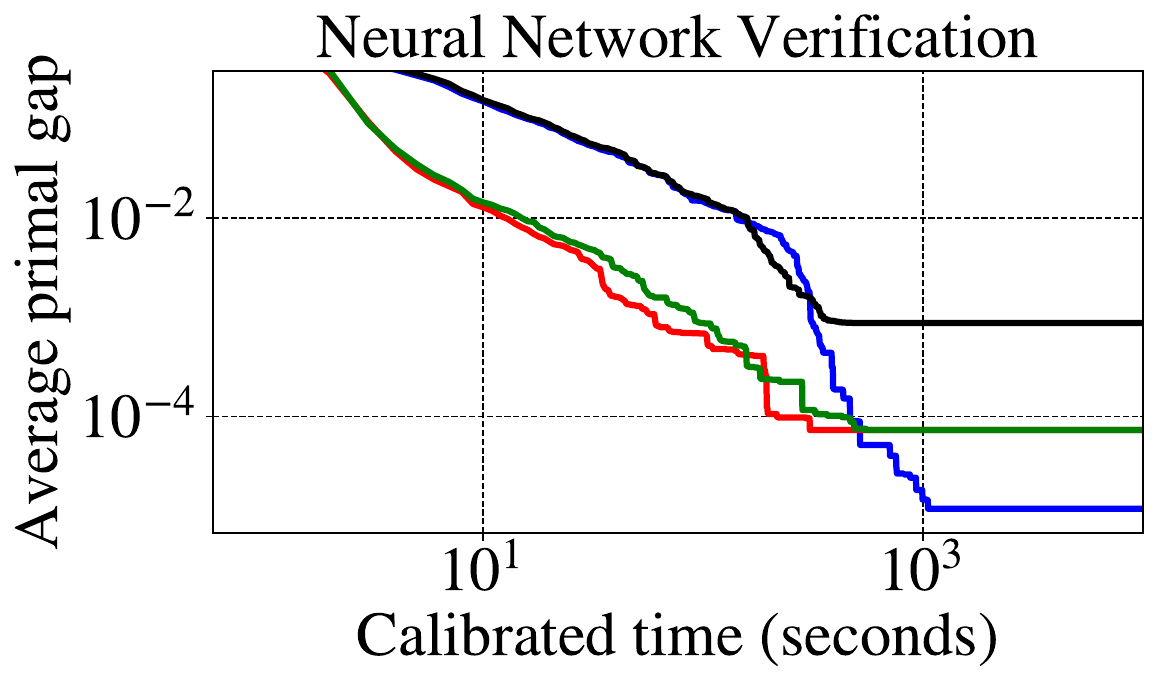}
        \includegraphics[width=\linewidth]{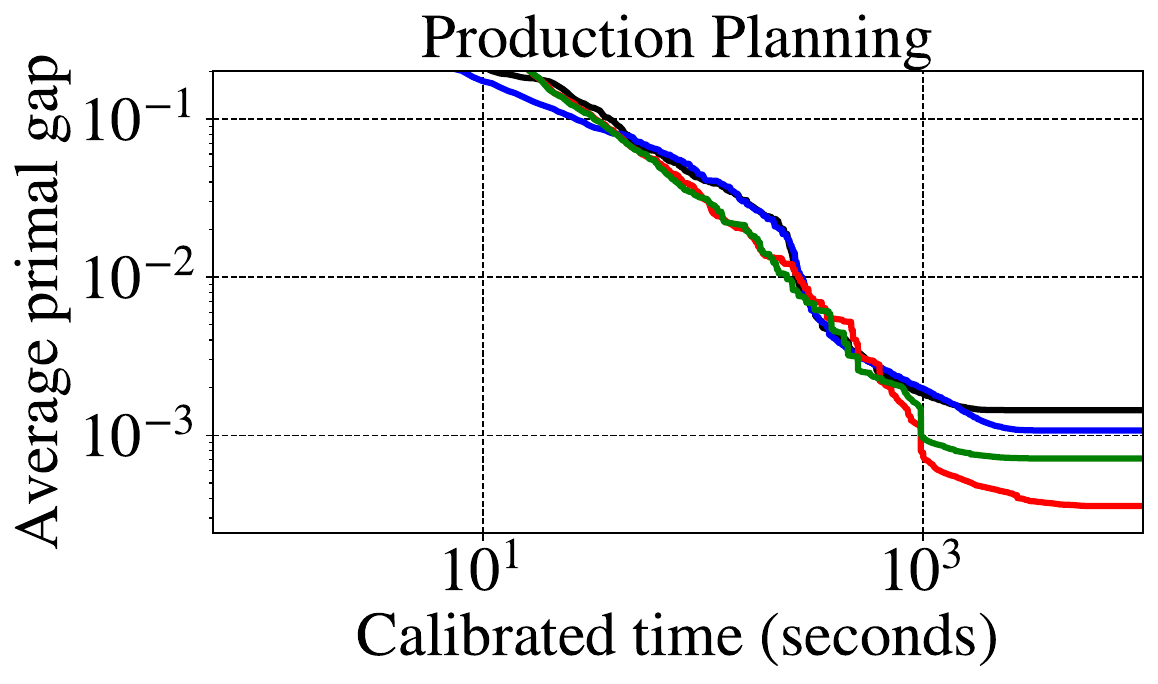}
        \includegraphics[width=\linewidth]{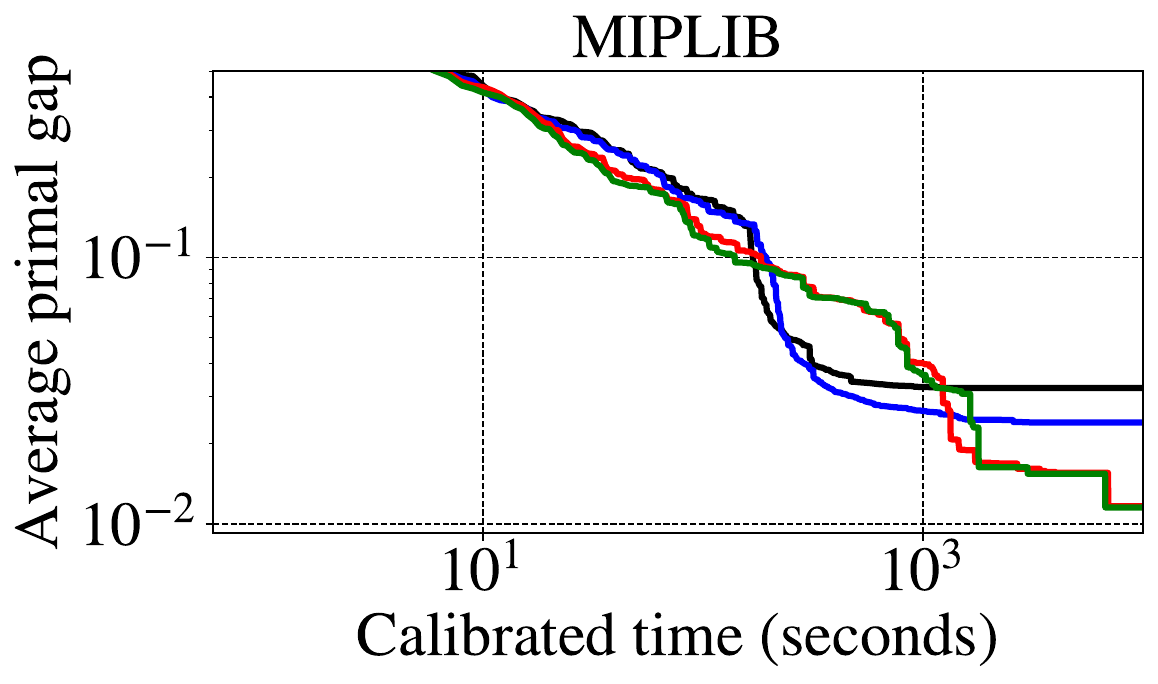} 
        \caption{Test set average primal gap as a function of running time for five datasets, for four combinations of two approaches to find an initial assignment (SCIP vs. Neural Diving), and two approaches to select a neighborhood (Neural vs. Random Neighborhood Selection).}
        \label{fig:gap_ablation}
    \end{subfigure}
    \label{fig:experiments}
    \caption{Experimental results.}
\end{figure*}

Figure~\ref{fig:survival} shows the performance of ND + NNS using survival curves. Compared to SCIP, our method's performance is considerably stronger on Production Packing, Electric Grid Optimization, and MIPLIB. On the first two, NNS solves almost all test instances to within the specified target gap, while SCIP only solves about 10\% on Production Packing, and about 80\% on Electric Grid Optimization. For Neural Network Verification, while SCIP eventually also solves all the instances, the survival curve for ND + NNS achieves the same fraction of solved instances faster. Even on MIPLIB, ND + NNS achieves a final solve fraction of roughly 80\%, compared to SCIP's 60\%. Similarly, comparing ND + NNS to Neural Diving shows the former achieving higher final solve fractions on all datasets except Production Planning, where the two methods perform roughly the same.

ND + NNS outperforms ND + RNS on Electric Grid Optimization, Neural Network Verification, and MIPLIB, by either achieving a better final solve fraction or the same solve fraction in less time. Note that survival curves need not fully reflect the improvements in average primal gaps achieved by ND + NNS shown in figure~\ref{fig:primal_gap} because improving the gap beyond the threshold does not improve the survival curve.

\begin{SCfigure*}
    \centering
    \begin{tabular}{cc}
        \includegraphics[width=0.3\textwidth]{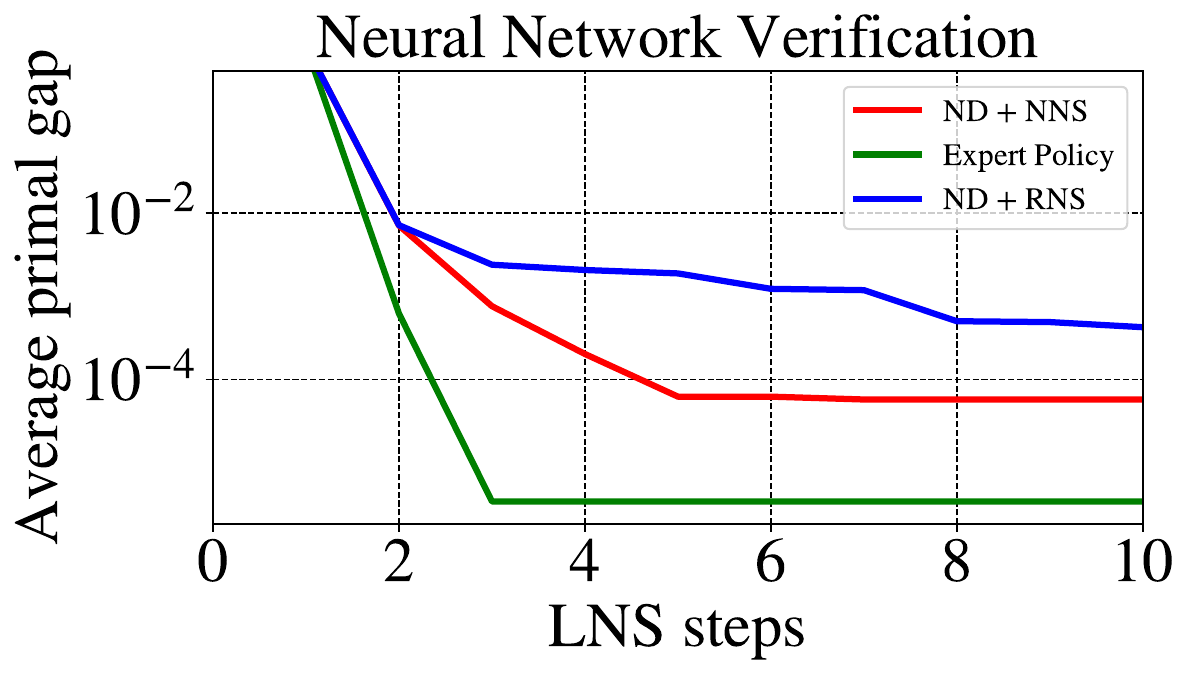} &
        \includegraphics[width=0.3\textwidth]{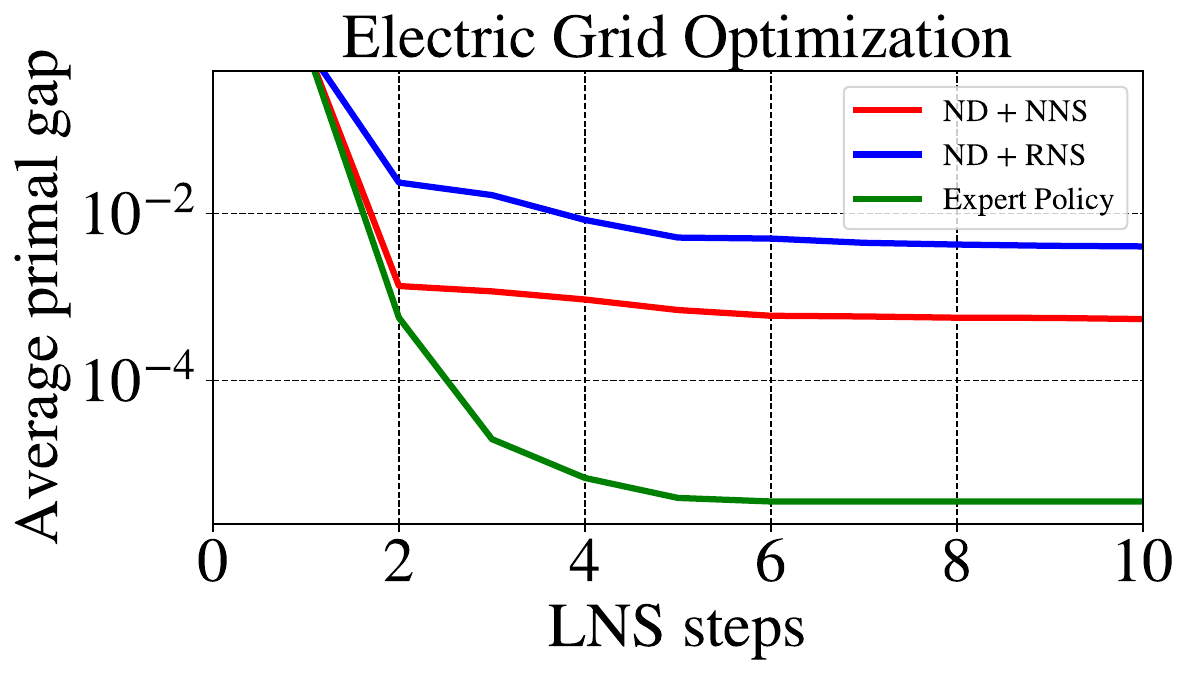} 
    \end{tabular}
    
    \caption{Comparison of expert policy used as the target for imitation learning to random (ND + RNS) and learned (ND + NNS) policies for selecting a search neighborhood at each step of large neighborhood search, with the initial assignment computed using Neural Diving for all three cases.}
    \label{fig:gap_by_step}
    
\end{SCfigure*}

\subsection{Ablation Study}
\label{subsec:ablation}

We evaluate how the two main components of our approach contribute to its performance. We consider four variants in which the initial assignment is given by either SCIP or Neural Diving, and neighborhood search is done using either Random Neighborhood Selection or Neural Neighborhood Selection. Figure~\ref{fig:gap_ablation} shows that, on all datasets except Neural Network Verification and MIPLIB, the average primal gap becomes worse without Neural Diving. This is true regardless of whether we use NNS or RNS. For MIPLIB, ND + NNS finishes with the best average primal gap, but is worse at intermediate running times. For Neural Network Verification, SCIP turns out to be better than Neural Diving for providing the initial assignment. NNS is crucial, achieving roughly a $100\times$ lower gap than SCIP + RNS. While the relative contribution of Neural Diving and Neural Neighborhood Selection to our approach's performance depends on the dataset, it is clear across all datasets that learning is necessary to achieve the best performance.

\subsection{Approximating the Expert}

Figure \ref{fig:gap_by_step} shows the average primal gap as a function of the number of LNS steps for the expert policy, Neural Neighborhood Selection, and Random Neighborhood Selection, all initialized with the same assignment computed using Neural Diving. This allows comparing the quality of the policies independently of their speed. For brevity we include only two representative datasets. The average primal gap is computed on a subset of the validation set. On both datasets the expert reduces the gap in the fewest steps, confirming its effectiveness as an imitation target. NNS learns to approximate it well enough to outperform RNS, which shows that imitation learning is a viable approach for the neighborhood selection task. Figure \ref{fig:gap_by_step} also shows that there is more room for the learned policy to approximate the expert better. Improving the learned policy's approximation while maintaining its speed can lead to even better average primal gap.

\section{Discussion}
\label{ref:discussion}
Neural Diving and ND + NNS have complementary strengths. As shown in section \ref{sec:results}, Neural Diving is often quicker than SCIP at achieving low average primal gaps in short running times, but its average primal gap tends to plateau at higher running times. One reason is that a single run of Neural Diving applies the learned model on the input MIP only once at the beginning, which has a fixed time cost. Any additional running time available is not used to apply the learned model again; instead it is allocated to the sub-MIP solve with SCIP. ND + NNS on the other hand can use the available time to repeatedly apply the NNS policy with more LNS steps. So higher running time limits can exploit the learned model better. As a result, ND + NNS is able to improve on Neural Diving at higher running time limits.

\noindent\textbf{Limitations:} 1) Our approach does not currently train models in an end-to-end fashion to directly optimize a final performance metric such as the average primal gap. Approaches based on Reinforcement Learning (RL) and offline RL may address this limitation. 2) The conditionally-independent model (section \ref{subsec:policy_network_architecture}) does not approximate the expert policy perfectly (figure \ref{fig:gap_by_step}). While this architecture supports efficient inference, the restrictive conditional independence assumption also limits its capacity. More powerful architectures, e.g., autoregressive models combined with GCNs, may give better performance by approximating the expert better.

\section{Related Work}
\label{sec:related_work}

\noindent\textbf{Learning and MIP Primal Heuristics:} \citet{ding2019accelerating} trains a neural network to predict values for a subset of a MIP's binary variables, which are then used to significantly reduce the search space by defining an additional constraint that any assignment be within a pre-specified Hamming distance from the predicted values. Similarly, \citet{xavier2020mluc} learns to warm start a MIP solver specifically for the electric grid \emph{unit commitment} problem using a $k$-Nearest Neighbor binary classifier to predict values for a subset of binary variables in the MIP. Both works are similar in spirit to Neural Diving, and shares its weakness that the model is applied in a one-shot manner at test time without the ability to iteratively improve the assignment. ND + NNS addresses this weakness.

Learning has also been used to choose among an ensemble of \emph{existing} primal heuristics in a complete solver. \citet{khalil2017primal} learn a binary classifier to predict whether applying a primal heuristic at a search tree node will improve the current best assignment. \cite{Hendel2018alns} formulate a multi-armed bandit approach to learn a switching policy online. This is complementary to our approach of constructing neural primal heuristics for a given application and can be combined by adding the neural heuristics to the ensemble.

\noindent\textbf{Other Learning Approaches for MIPs:} Several works have used learning to improve tree search in a complete solver \cite{He2014LearningToSearch, Khalil2016LearningToBranch, Alvarez2017MLApproxStrongBranching, gasse2019exact, zarpellon2020parameterizing, gupta2020hybrid}. They are most relevant when both an assignment and its optimality gap are required. Here we consider the setting where only the former is required. Learning has been used to predict hyperparameters for MIP solvers, either for \emph{Algorithm Configuration} \citep{ansotegui2009gga, hutter2009paramils, hutter2011smac, ansotegui2015ggapp} or \emph{Algorithm Selection} \citep{kotthoff2016AlgoSelectSurvey, hutter2014epm}. Such approaches can be combined with ours.

\noindent\textbf{Learning for Combinatorial Optimization:} Our work is an instance of learning for combinatorial optimization problems. Some of the earliest works in this area are \cite{zhang1995RLforJSS, moll1999routing, boyan97stage}. More recently, deep learning has been applied to the Travelling Salesman Problem \citep{vinyals2015pointernets, bello2016neural}, Vehicle Routing \citep{kool2018attention, nazari2018RLforVRP}, Boolean Satisfiability \citep{selsam2019NeuroSAT, amizadeh2018learning, yolcu2019localsearch}, and general graph-structured combinatorial optimization problems \citep{khalil2017learning, li2018combinatorial}. A survey of the topic is available by \citet{bengio2018survey}.
\section{Summary \& Conclusions}
\label{ref:conclusions}
We have proposed a learning-based LNS approach for MIPs. It trains a Neural Diving model to generate an initial assignment, and a Neural Neighborhood Selection policy to select a search neighborhood at each LNS step. The resulting neighborhood can be searched by solving with SCIP a smaller sub-MIP derived from the input MIP. Our approach matches or significantly outperforms all baselines with respect to average primal gap and survival curves on five datasets containing diverse, large-scale, real-world MIPs. It addresses a key limitation of Neural Diving as a standalone primal heuristic by improving the average primal gap at larger running times. Even larger performance gains can potentially be achieved with end-to-end training of both models, and using more powerful network architectures.

\bibliography{references}
\bibliographystyle{icml2022}

\newpage
\appendix
\onecolumn
\newpage

\section{Technical Appendix}
\label{sec:appendix}

\subsection{Expert data generation}
\label{subsec:expert_data_generation}
To generate trajectories of the expert policy, we computed 10 steps of the local branching procedure as described in Section 3.2 for each instance in the training datasets. The Hamming distance $\eta$ of the local branching step was picked as a fixed fraction of the number of variables $n$, i.e. $\eta=r\cdot n$. For each step computation we used SCIP with default parameters, with a time limit of 3 hours. The parameter $r$ was chosen using the validation set, and the results are shown in table \ref{tab:train_params}.

\subsection{Training procedure}
\label{subsec:training_procedure}
A model was trained for each dataset on a single NVIDIA V100 GPU for up to 400k steps, taking up to 48 hours. We used grid-search to tune hyperparameters (number of GCN layers $n_{l}$, learning rate (lr) and its decay, where we decayed the learning rate by 0.9 every $m$ steps) on the validation set. The resulting choices, as well as the training set size $N$ are shown in the table below. Node embeddings were of size $64$ for each model, and were computed using a $2$-layer MLP with $64$ hidden units per layer.

\begin{table}
\centering
\begin{tabular}{c|c|c|c|c|c}
    \bf Dataset & $N$ & $r$ & $n_l$ & lr & $m$ \\
    \hline
    Prod. Planning & 7189 & 0.05 & 5 & $5 \times 10^{-4}$ & 500 \\
    Prod. Packing & 2924 & 0.01 & 20 & $10 ^ {-3}$ & 500 \\
    MIPLIB & 553 & 0.05 & 5 & $10^{-3}$ & 500 \\
    NN Verification & 2554 & 0.4 & 7 & $10^{-3}$ & 100 \\
    Electric Grid Opt. & 11444 & 0.4 & 10 & $10^{-3}$ & 500 \\
\end{tabular}
\caption{Hyperparameters for data generation and training.}\label{tab:train_params}
\end{table}

\subsection{Evaluation procedure}
\label{subsec:evaluation_procedure}
During evaluation, we used the trained model as described in Section 3.5. For all datasets, we performed 50 steps of LNS, and the adaptation factor for the neighbourhood size was fixed to $\alpha=1.5$ and sampling probability bias to $\epsilon=10^{-3}$.
Hyperparameters (initial fraction of unassigned variables $r_0$, sampling temperature $\tau$, and search time limit per step $t$) were optimized using the validation set to the values in table \ref{tab:eval_params}.

\begin{table}
\centering
\begin{tabular}{c|c|c|c}
    \bf Dataset & $r_0$ & $\tau$ & \bf  $t$ (sec) \\
    \hline
    Prod. Planning & 0.05 & 2.0 & 300 \\
    Prod. Packing & 0.15 & 1.0 & 300 \\
    MIPLIB & 0.25 & 2.0 & 300 \\
    NN Verification & 0.4 & 1.5 & 60 \\
    Electric Grid Opt. & 0.4 & 1.5 & 300 \\
\end{tabular}
\caption{Hyperparameters for evaluation.}\label{tab:eval_params}
\end{table}

\subsection{Input representation}
\label{subsec:input_representation}
We use input features from \cite{gasse2019exact} along with the incumbent solution. The code for computing the features from \cite{gasse2019exact} is available at \href{https://github.com/ds4dm/learn2branch}{https://github.com/ds4dm/learn2branch}.

\subsection{Calibrated time}
\label{appendix:calibrated_time}

The total evaluation workload across all datasets and comparisons requires a large amount of compute. To meet the compute requirements, we use a shared, heterogeneous compute cluster. Accurate running time measurement on such a cluster is difficult because the tasks may be scheduled on machines with different hardware, and interference from other unrelated tasks on the same machine increases the variance of solve times. To improve accuracy, for each solve task, we periodically solve a small \emph{calibration MIP} on a different thread from the solve task on the same machine. We use an estimate of the number of calibration MIP solves during the solve task on the same machine to measure time, which is significantly less sensitive to hardware heterogeneity and interference. This quantity is then converted into a \emph{calibrated time} value using the calibration MIP's solve time on a reference machine.

We define the speed of a machine on which a MIP solving job runs to be
\begin{equation}
    \text{Speed} = \frac{1}{\text{Wall clock time to solve calibration MIP}}.
\end{equation}
For each periodic measurement of calibrated time, we estimate the speed $K$ times and use the average. $K$ is set to be the number of samples needed to estimate mean speed with 95\% confidence, with a minimum of 3 samples and a maximum of 30. The elapsed calibrated time $\Delta t_{\text{calibrated}}$ since the last measurement is
\begin{equation}
\label{eqn:calibrated_time}
    \Delta t_{\text{calibrated}} = \text{Speed} \times \Delta t_{\text{wallclock}},
\end{equation}
where $\Delta t_{\text{wallclock}}$ is the elapsed wallclock time since the last measurement. We use the MIP named \emph{vpm2} from MIPLIB2003 \cite{miplib2003} as the calibration MIP.

Note that the above definition of calibrated time does not have a time unit. Instead it is (in effect) a count of the calibrated MIP solves during the evaluation solve task. To give it a unit of seconds, one can choose a reference machine with respect to which evaluation solve times will be reported, accurately measure the calibration MIP's solve time on it (without other tasks interfering), and multiply the calibrated time in equation \ref{eqn:calibrated_time} by the reference machine's estimated calibration MIP solve time. The resulting quantity has a unit of seconds. It can be interpreted as the time the evaluation solve task would have taken if it ran on the reference machine. We select Intel Xeon 3.50GHz CPU with 32GB RAM as the reference machine. All the calibrated time results in the paper are expressed with respect to the reference machine in seconds.

\subsection{Dataset Details}
\label{subsec:appendix_datasets}

The target optimality gap used for each dataset in our evaluation as SCIP's stopping criterion for a MIP solve is given in table~\ref{tab:target_opt_gaps}. In the case of Neural Network Verification, the actual criterion used in the application is to stop when the objective value becomes negative, but this is not expressible as a constant target gap across all instances. In order to treat all datasets consistently, we have selected a gap of 0.05. Additional information about the applications that the datasets are extracted from is provided in table~\ref{table:datasets}. 

\begin{table}
    \caption{Optimality gap thresholds used for plotting survival curves for the datasets in our evaluation.}
    \centering
    \begin{tabular}{c|c}
        \toprule
        Dataset & Target Optimality Gap \\
        \midrule
        Neural Network Verification & 0.05\\
        Production Packing & 0.01 \\
        Production Planning & 0.03\\
        Electric Grid Optimization & 0.0001 \\
        MIPLIB & 0\\
        \bottomrule
    \end{tabular}
    \label{tab:target_opt_gaps}
\end{table}

\begin{figure*}
\centering
    \includegraphics[width=0.46\textwidth]{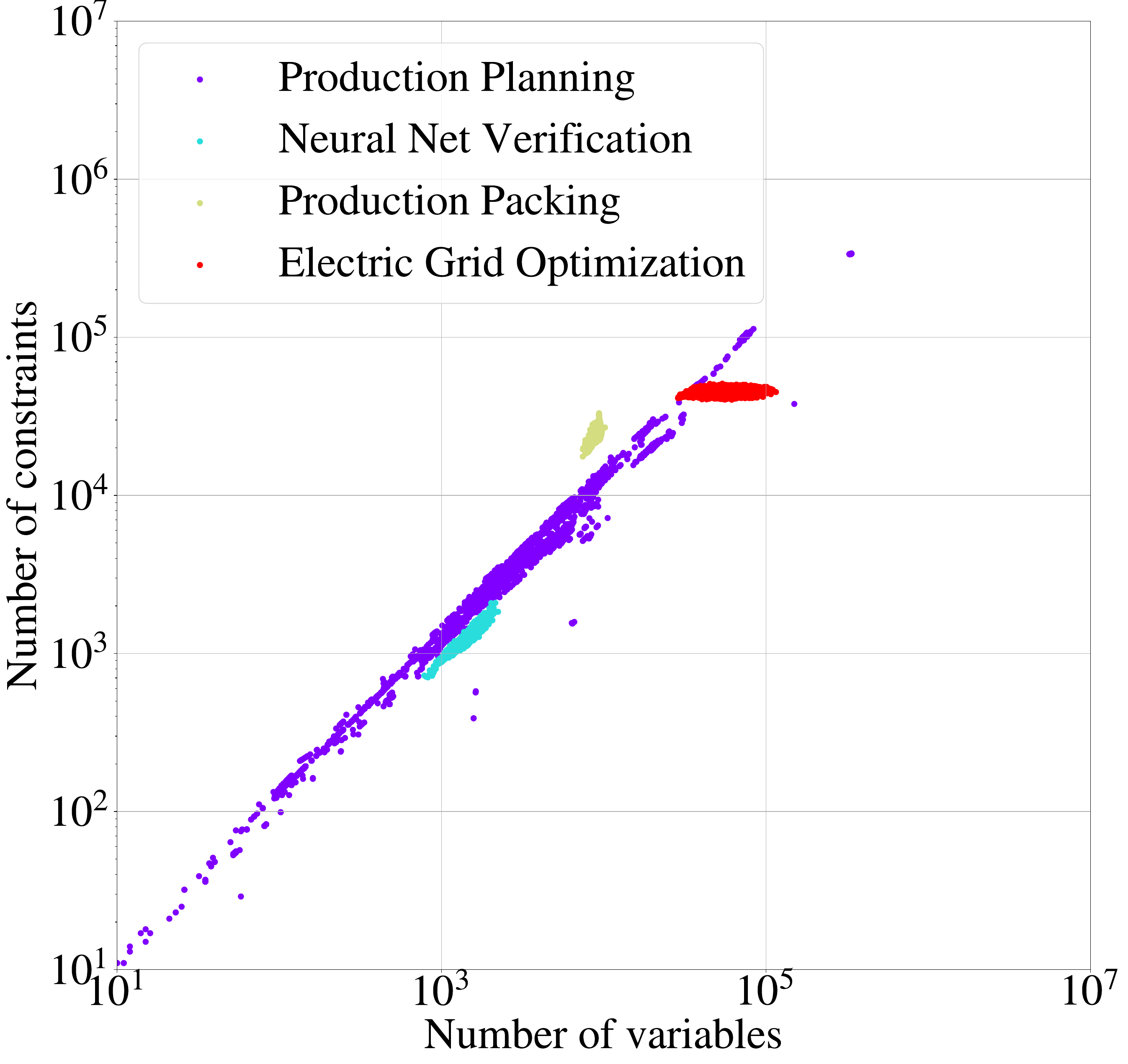}
    \hspace{0.1in}
    \includegraphics[width=0.46\textwidth]{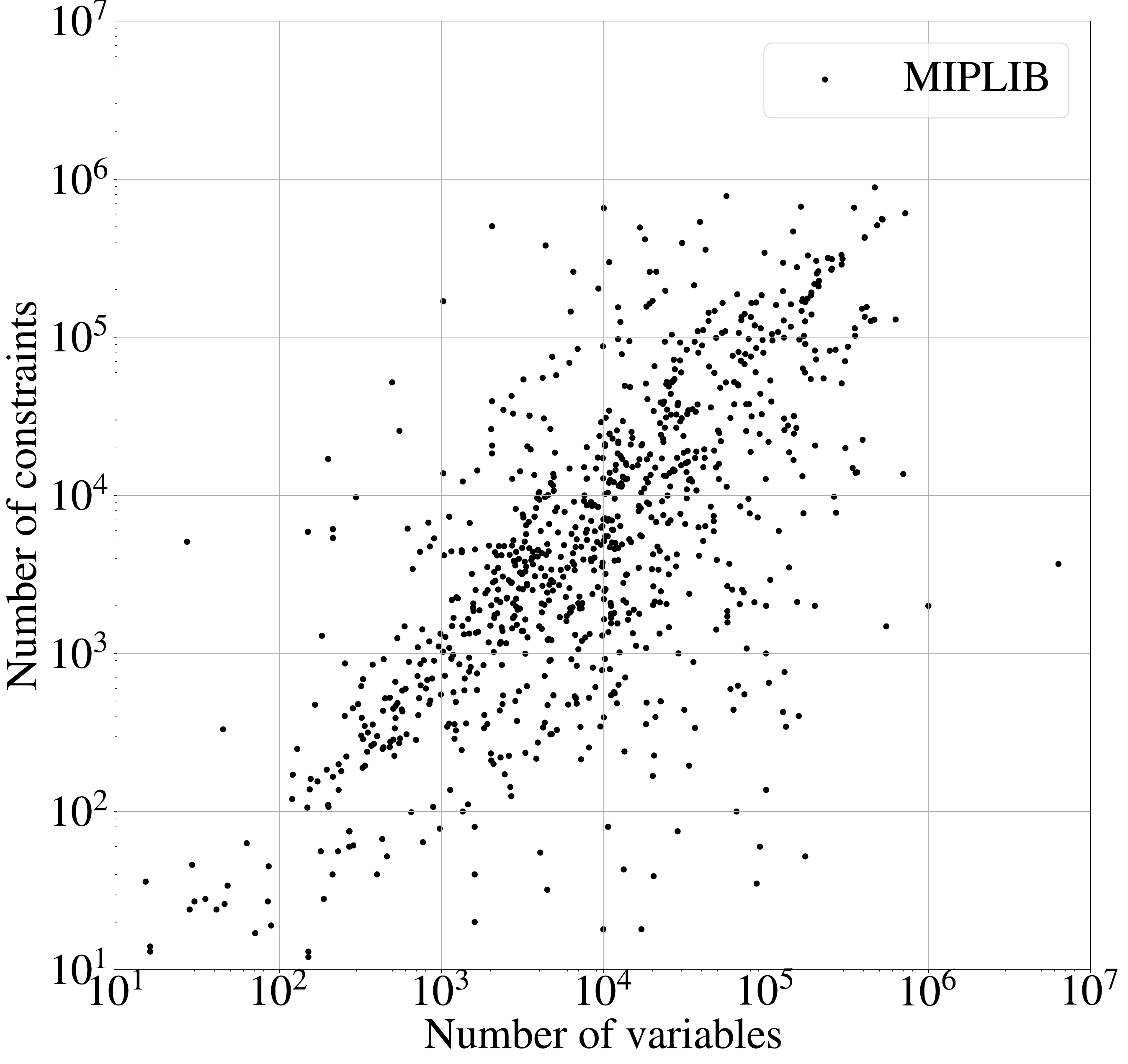}
\caption{
    {\small Number of variables versus number of constraints \emph{after presolving} using SCIP 7.0.1 for the application-specific datasets (left) and MIPLIB (right) used in our evaluation. Presolving significantly reduces the problem size compared to that of the raw input MIP.}
}
\label{fig:presolved_mip_sizes}
\end{figure*}

Figure \ref{fig:presolved_mip_sizes} shows the MIP sizes for the datasets used in our evaluation after presolving using SCIP 7.0.1. Note that, among the application-specific datasets, the Production Planning dataset is the most heterogeneous (along with MIPLIB) in terms of instance sizes. Table~\ref{tab:dataset_stats} summarizes the characteristics of those datasets before and after presolving with the SCIP 7.0.1 solver. The number of constraints and variables ranges different orders of magnitude across the datasets. It is worth noting that during presolving some instances might be deemed infeasible and, hence, dropped from the dataset.

\begin{table*}[t]
\caption{Description of the five datasets we use in the paper. Please see \cite{nair2020solving} for more details.}
\label{table:datasets}
\begin{tabular}{|p{4cm}|p{13cm}|}
\hline
    \textbf{Name} & \textbf{Description} \\
    
    \hline
    Neural Network Verification & Verifying whether a neural network is robust
    to input perturbations can be posed as a MIP \citep{cheng2017verification, tjeng2018evaluating}.
    Each input on which to verify the network gives rise to a different MIP.
    In this dataset, a convolutional neural network is verified on each
    image in the MNIST dataset, giving rise to a corresponding dataset of MIPs.\\ %

    \hline

    Production Packing & A packing optimization problem solved in a large-scale production system.\\
    
    \hline
    
    Production Planning & A planning optimization problem solved in a large-scale production system.\\
    
    \hline
    
    Electric Grid Optimization & Electric grid operators optimize the choice of power generators to use at different times of the day to meet electricity demand by solving a MIP. This dataset is constructed for one of the largest grid operators in the US, PJM, using publicly available data about generators and demand, and the MIP formulation in \citep{knueven2018mip_uc}.\\
    
    \hline
    
    MIPLIB & Heterogeneous dataset containing `hard' instances of MIPs
    across many different application areas that is used as a long-standing standard benchmark for MIP solvers \citep{gleixner2019miplib}. We use instances from both the 2010 and 2017 versions of MIPLIB. \\ 
\hline
\end{tabular}
\end{table*}

\begin{table*}
    \caption{Statistics (\textbf{median} in the first row and \textbf{maximum} in the second row) of constraints and variables (per type) for the different datasets before / \textit{after} presolving using SCIP 7.0.1. The total number of nodes in the bipartite graph representation corresponds to the sum of the number of constraints and variables.}
    \begin{adjustbox}{center}
    {\footnotesize
    \begin{tabular}{ccccccc}
    \toprule
        Dataset & Constraints & Variables & Binary & Integer (non-binary) & Continuous  \\
        \midrule 
        \emph{Neural Network Verification} & 6531/\textit{1407} & 7142/\textit{1629} & 171/\textit{170} & 0/\textit{0} & 6972/\textit{1455} \\
        (max) & 7123/\textit{2089} & 7535/\textit{2231} & 364/\textit{364} & 0/\textit{0} & 7171/\textit{1921} \\
         \midrule
        \emph{Production Packing} & 36905/\textit{24495} & 10046/\textit{8919} & 3773/\textit{3437} & 0/\textit{0} & 6231/\textit{5421} \\
        (max) & 47414/\textit{33071} & 11233/\textit{10161} & 4120/\textit{3771} & 0/\textit{0} & 7113/\textit{6390} \\
        \midrule
        \emph{Production Planning} & 11910/\textit{478} & 9884/\textit{404} & 833/\textit{119} & 462/\textit{119} & 8337/\textit{136} \\
        (max) & 1722678/\textit{338508} & 1548582/\textit{335783} & 117485/\textit{22626} & 58782/\textit{23628} & 1374182/\textit{335783} \\
         \midrule
        \emph{Electric Grid Optimization} & 61851/\textit{45834} & 60720/\textit{58389} & 42240/\textit{42147} & 0/\textit{0} & 18768/\textit{16623} \\
         (max) & 67086/\textit{50908} & 117792/\textit{115098} & 98112/\textit{98021} & 0/\textit{0} & 21456/\textit{20184} \\
        \midrule
        \emph{MIPLIB} & 7706/\textit{4388} & 11090/\textit{9629} & 4450/\textit{3816} & 0/\textit{0} & 218/\textit{96} \\
         (max) & 19912111/\textit{888363} & 38868107/\textit{6338552} & 20677405/\textit{6338552} & 549428/\textit{0} & 38335410/\textit{700426} \\
        \bottomrule
    \end{tabular}
    }
    \end{adjustbox}
    \label{tab:dataset_stats}
\end{table*}

\subsection{Results on NeurIPS'21 Competition Datasets}
\label{subsec:neurips21_results}
\begin{figure*}
    \centering
    \begin{subfigure}[t]{0.45\textwidth}
        \centering
        \includegraphics[width=\linewidth]{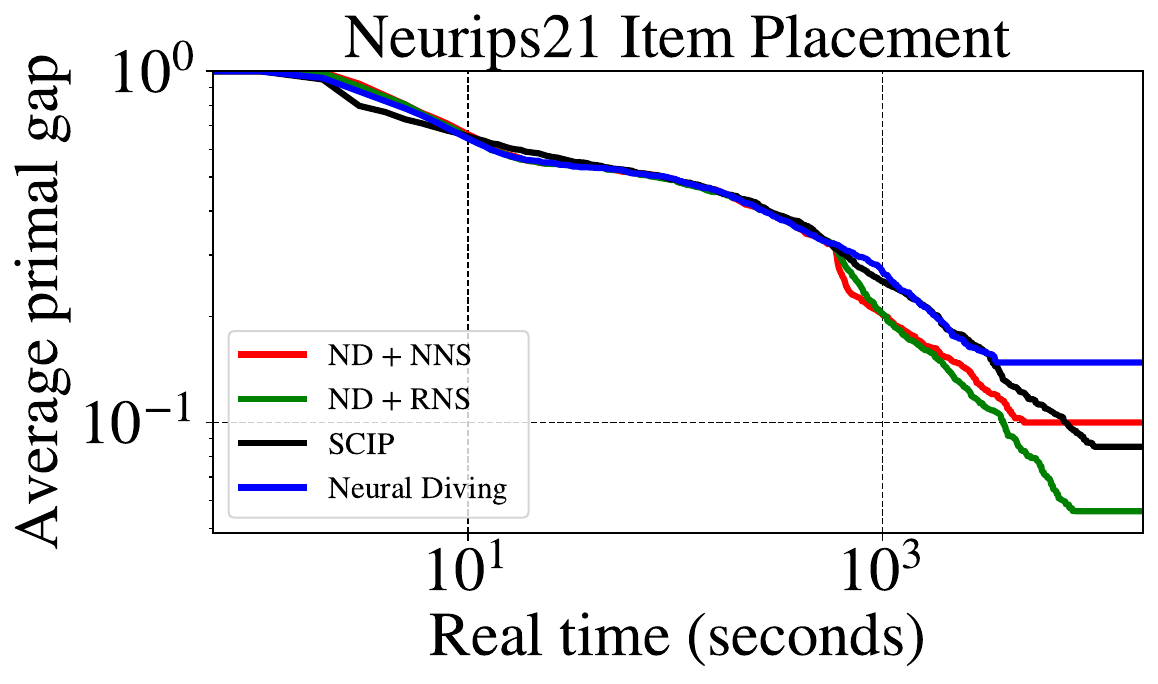}
        \includegraphics[width=\linewidth]{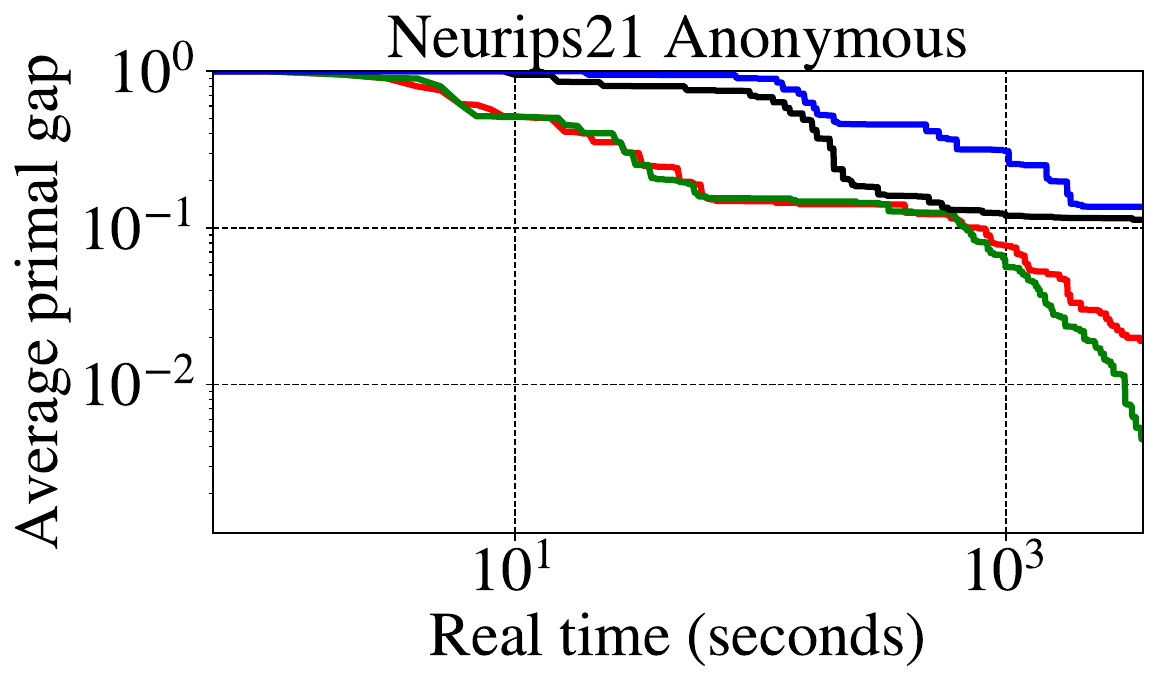}
        \includegraphics[width=\linewidth]{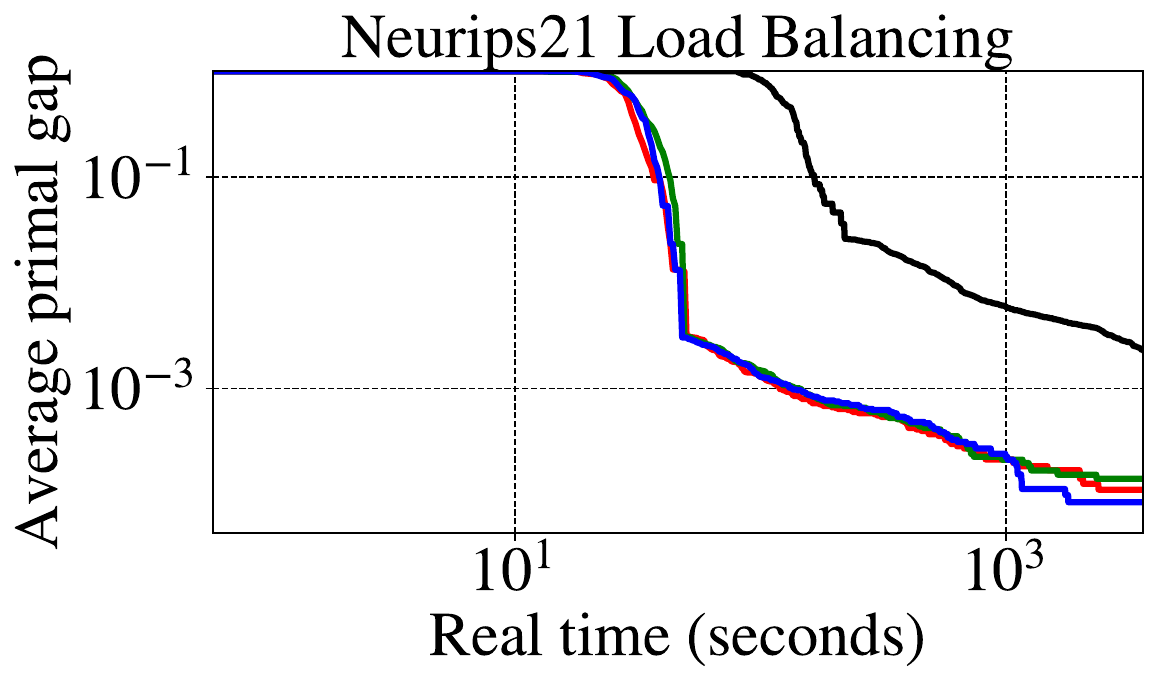}        
        \caption{Fraction of test set instances with primal gap below a dataset-specific threshold, as a function of running time for five datasets. (Note: For Production Planning, several curves closely overlap.)}
        \label{fig:primal_gap_neurips21}
    \end{subfigure}
    \hfill
    \begin{subfigure}[t]{0.45\textwidth}
        \centering
        \includegraphics[width=\linewidth]{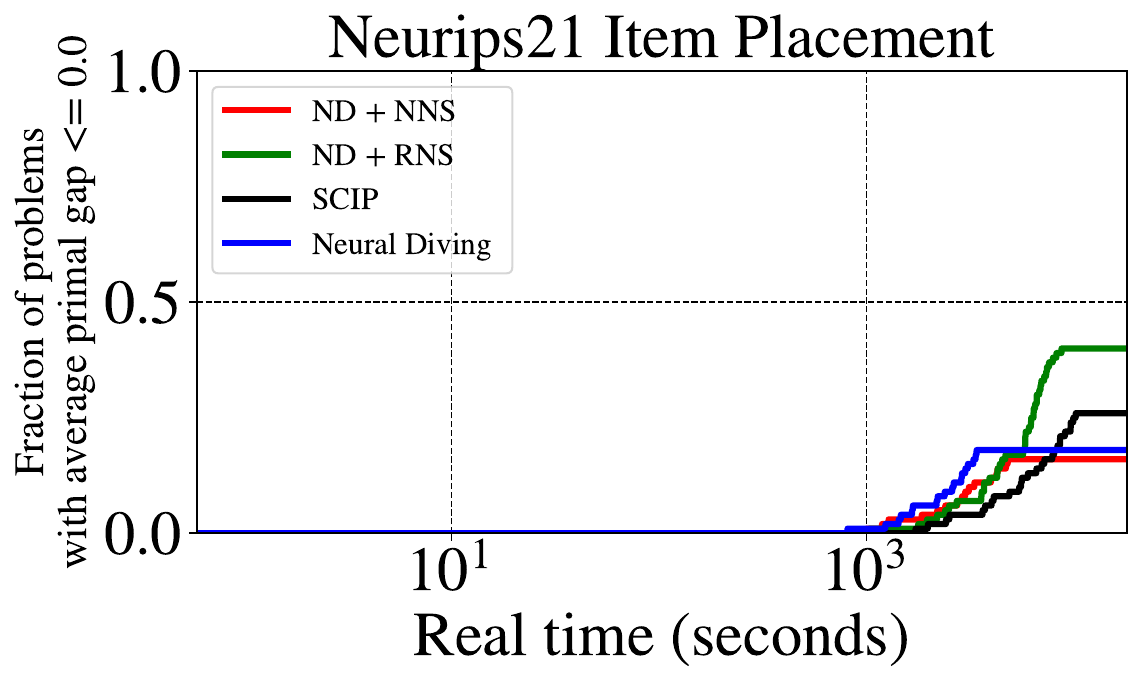} 
        \includegraphics[width=\linewidth]{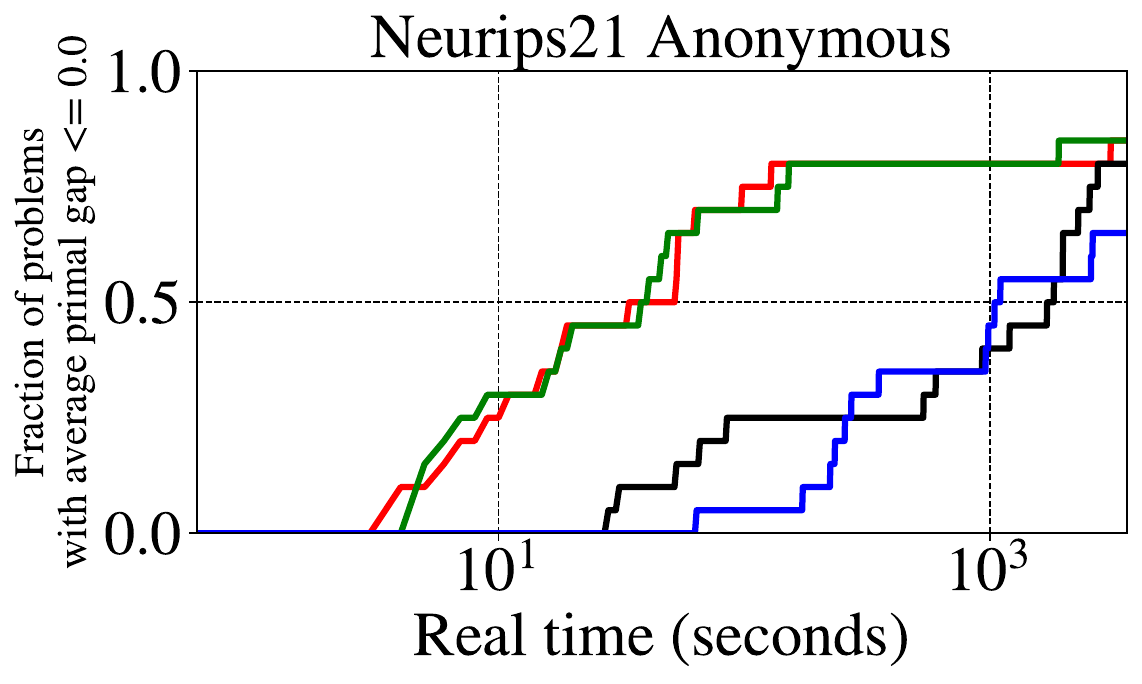} 
        \includegraphics[width=\linewidth]{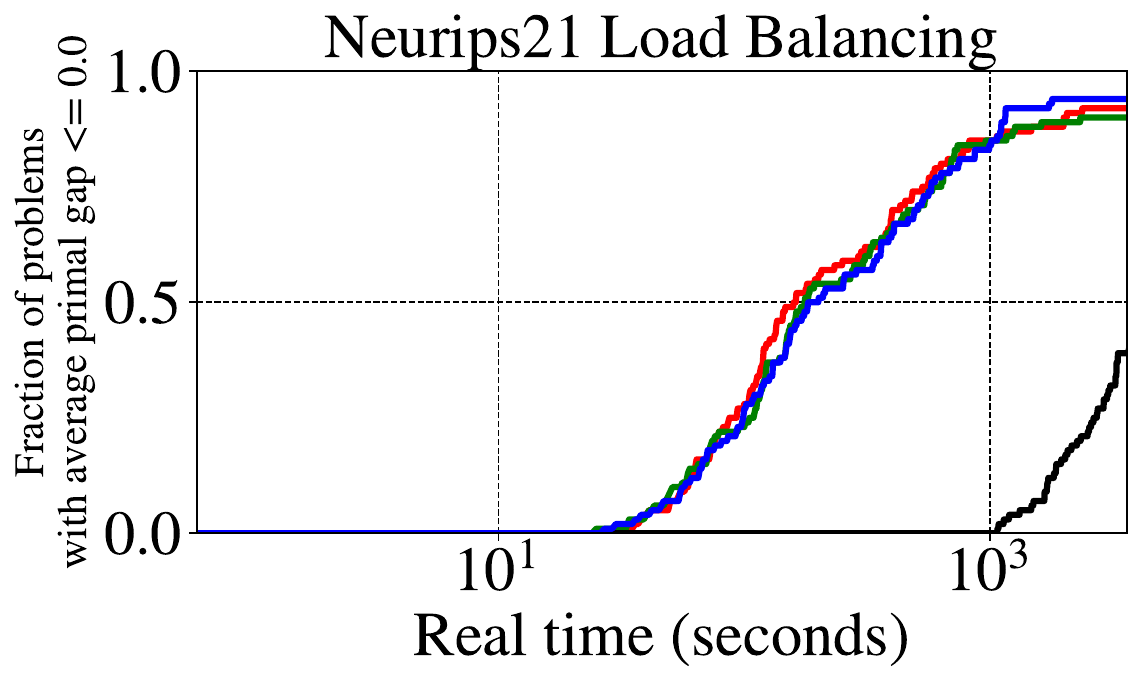}      
        \caption{Test set average primal gap (see section \ref{subsec:metrics}, lower is better) as a function of running time for five datasets.}
        \label{fig:survival_neurips21}
    \end{subfigure}
    \caption{Results for the NeurIPS'21 Machine Learning for Combinatorial Optimization competition datasets.}
    \label{fig:neurips21_experiments}    
\end{figure*}

Figure \ref{fig:neurips21_experiments} presents preliminary results based on limited hyperparameter tuning for Neural Diving and Neural Neighborhood Search on the three datasets (Item Placement, Anonymous, Load Balancing) from the NeurIPS'21 competition on \href{https://www.ecole.ai/2021/ml4co-competition/}{Machine Learning for Combinatorial Optimization}. The Anonymous dataset is very small, with only 98 training cases. Not surprisingly, Neural Diving does not perform well on this dataset, and Neural Neighborhood Search does not provide any benefit over Random Neighborhood Search. Item Placement is a challenging dataset for SCIP. We have observed that even with several hours of running time for SCIP during data collection, SCIP still frequently reports very high optimality gaps (often $> 100\%$). This affects both the quality of the data collected to train the ND and NNS models, as well as the test time performance of ND and NNS by using SCIP to solve the sub-MIPs, which explains the poor performance of all approaches on this dataset.  On Load Balancing, ND significantly outperforms SCIP, which shows that learning is useful for this dataset. However, NNS and RNS both fail to provide any improvements over ND. We believe further tuning the hyperparameters can improve the results on Load Balancing.

\end{document}